\documentclass[english,preprint,11pt]{article}
\usepackage[latin9]{inputenc}
\usepackage[T1]{fontenc}
\usepackage{calc}
\usepackage{url}
\usepackage{amsmath}
\usepackage{amssymb}
\usepackage{graphicx}
\usepackage{setspace}
\usepackage{wrapfig}
\usepackage{natbib}
\usepackage{color}
\usepackage{comment}
\usepackage{array}
\usepackage{multirow}
\usepackage{pdflscape}
\usepackage[normalem]{ulem}
\usepackage[affil-it]{authblk}
\usepackage[top=1in, left=1in, right=1in, bottom=1in]{geometry}
\makeatletter
\newcommand\textlcsc[1]{\textsc{\MakeLowercase{#1}}}

\usepackage{babel}

\providecommand{\tabularnewline}{\\}

\DeclareMathOperator*{\argmin}{\mbox{argmin}}
\DeclareMathOperator*{\argmax}{\mbox{argmax}}
\newenvironment{hanglist} 
{\begin{list}{}{\setlength{\itemsep}{4pt}
			\parindent 0pt  \setlength{\parsep}{0pt}\setlength{\leftmargin}{+25pt}
			\setlength{\itemindent}{-\parindent}}}{\end{list}}

\begin{document}

\title{Branch-and-price algorithm for an auto-carrier transportation problem}

\author{Saravanan Venkatachalam%
	\thanks{email: \texttt{saravanan@tamu.edu}; Corresponding author}}
\affil{Engineering Technology and Industrial Distributuion Department, Texas A\&M University,\\ College Station, Texas 77843}

\author{Kaarthik Sundar%
	\thanks{email: \texttt{kaarthiksundar@tamu.edu}}}
\affil{Department of Mechanical Engineering, Texas A\&M University, \\ College Station, Texas 77843}

\date{Dated: \today}

\maketitle

\begin{abstract}
Original equipment manufacturers (OEMs) manufacture, inventory and transport new vehicles to franchised dealers. These franchised dealers inventory and sell new vehicles to end users. OEMs rely on logistics companies with a special type of truck called an auto-carrier to transport the vehicles to the dealers. The process of vehicle distribution has a common challenge. This challenge involves determining routes, and the way to load the vehicles onto each auto-carrier. In this paper, we present a heuristic to determine the route for each auto-carrier based on the dealers' locations, and subsequently, a branch-and-price algorithm to obtain optimal solutions to the loading problem based on the generated route. The loading problem considers the actual dimensions of the vehicles, and the restrictions imposed by vehicle manufacturers and governmental agencies on the loading process. We perform extensive computational experiments for the loading problem using real-world instances, and our results are benchmarked with a holistic model to corroborate the effectiveness of the proposed method. For the largest instance comprising of 600 vehicles, the proposed method computes an optimal solution for the loading problem within a stipulated runtime. \\


\end{abstract}

\noindent \hrulefill

\section{Introduction\label{sec:Introduction}}
The vehicle distribution system in the U.S. has a single dominant form, in which the original equipment manufacturers (OEMs) manufacture, inventory and transport new vehicles to franchised dealers. These franchised dealers inventory and sell new vehicles to the end customers \cite{Karabakal2000}. OEMs rely on logistics companies to distribute the vehicles to the dealers. Logistics companies have one or more distribution centers, and the manufacturers transport the vehicles to these distribution centers. Apart being a transit point between the OEMs and dealers, the distribution centers process the vehicles for their minor customized requests (upgrade of seat covers, lights etc.) from the end customers. Distribution centers are ideal to perform such requests, otherwise the dealers will be burdened by inventory overhead or the OEMs' economy of scale will get affected. The logistics companies then deliver the vehicles to the dealers who in turn sell them to the end customers. This last mile of vehicle delivery to the dealers' locations is highly expensive. Taking into account that 16.4 million new vehicles were sold in the U.S. \citep{nada} for the year 2014, a conservative estimate of \$100 per vehicle for the last mile delivery puts the expected expenditure in excess of \$16 billion. Furthermore, the American Trucking Association accords this industry with a special status in their group \citep{acc}. 

\begin{figure}
	\centering
	\includegraphics[scale=0.65]{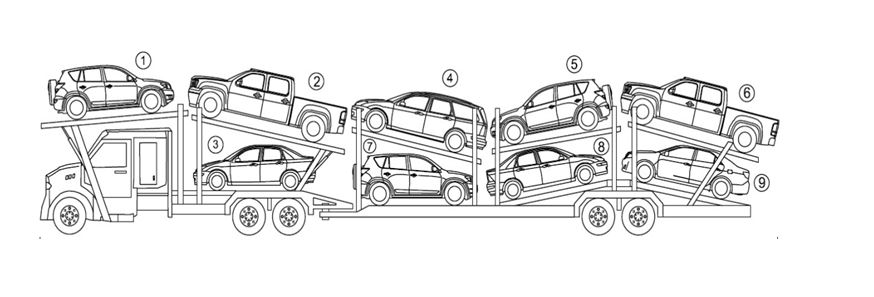}
	\caption{An auto-carrier with nine loading ramps.}
	\label{fig:illustration}
\end{figure}

This study presents a methodology to distribute the vehicles for an auto logistics company (ALC) owning auto-carriers. The proposed methodology is evaluated on a real case study of an ALC in the southern U.S. In this paper we denote \textit{vehicle} as an OEM's product (e.g., a car, a truck, or a SUV), \textit{auto-carrier} as a special type of truck owned by ALCs which is used to transport the \textit{vehicles}, and a \textit{dealer's location} is a dealership address for an OEM. The vehicles are transported from ALC's distribution center to the dealers' locations using  auto-carriers, and the delivery schedules to the dealers' locations are generated based on weekly demand. ALC serves dealers over a wide geographical area and to efficiently schedule the deliveries, the dealers are partitioned into clusters based on their locations. Each cluster may be further partitioned into sub-clusters based on the demand. We shall refer to these sub-clusters as regions. Instead of scheduling the deliveries to all its dealers, ALC schedules deliveries to the dealers in each region, and uses a heterogeneous fleet of auto-carriers for the deliveries. The dimensions and weight of the vehicles restrict the number of vehicles that can be loaded on each type of auto-carrier. A \textit{type} refers to a particular variety of auto-carriers as their ability to carry a certain number of vehicles differs from one type to another. An auto-carrier consists of a tractor and a trailer equipped with upper and lower loading ramps as shown in Figure \ref{fig:illustration}. Figure \ref{fig:illustration} illustrates an auto-carrier with nine ramps, which is typical. We will refer to the set of vehicles in a trip as \textit{load}. A load is always associated with an auto-carrier and has a precise vehicle-ramp assignment. Very few dealers in a given region have enough daily volume to receive an entire load, so a load may contain a mix of vehicles to be delivered to a group of dealers. The planner working at ALC is responsible for determining the routes and loads for the auto-carriers. We will refer to the planner as a \textit{user}. The primary objective of the user is to build loads to satisfy the demands of the dealers in a given region at a minimum cost. Given a set of vehicles to be delivered for a given region, the cost of the operations depends upon the number of auto-carriers used for the delivery. Hence, the planner's objective is to minimize the number of auto-carriers to deliver a given set of vehicles to the dealers in a region.
 
Auto-carriers are equipped with special loading equipment on the ramps to load as many vehicles as possible. For instance, in Figure \ref{fig:illustration}, loading ramps 6 and 9 can be extended horizontally, and ramps 2 and 9 can be rotated. Vehicles are loaded and unloaded on the auto-carriers from the rear, and unloading without reshuffling is preferable. The load of an auto-carrier depends on a number of constraints. The constraints include the restrictions on maximum length, height, and weight of the cargo set by the government authorities (U.S. Department of Transportation) \citep{dot}, and the number of reloads allowed. An instance of reloading would occur when a vehicle to be delivered is on an interior ramp of the auto-carrier, and the only way to remove the vehicle is by removing the vehicles on the outer ramps. Reloading increases the chances of accidental damage that can occur to a vehicle as it is being unloaded and reloaded, and is a time consuming process affecting the productivity of the resources. The authors of \cite{Agbegha1998} estimated the reloading cost to be in excess of \$22.5 million per year. To circumvent reloading, \cite{Agbegha1998} and  \cite{Dell2014} impose a Last-In-First-Out (LIFO) policy when loading the vehicles onto the ramps. On the downside, the LIFO policy tends to increase the number of trips for the logistics companies. Furthermore, a load can contain vehicles from different dealers, and if the sequence of deliveries is fixed a priori, then this affects the vehicle-ramp assignment of the load. 

\begin{table}[h]
	\begin{centering}
		{\small{}}%
		\begin{tabular}{|>{\centering}p{1.3cm}|c|c|c|c|c|c|c|c|c|}
			\hline 
			\multirow{3}{1.5cm}{} & \multicolumn{3}{c|}{{\small{Honda}}} & \multicolumn{3}{c|}{{\small{Toyota}}} & \multicolumn{3}{c|}{{\small{Ford}}}\tabularnewline
			\cline{2-10} 
			& {\small{Ridgeline}} & {\small{Accord}} & {\small{Fit}} & {\small{Tundra}} & {\small{Camry}} & {\small{Yaris}} & {\small{F350}} & {\small{Focus}} & {\small{Fiesta}}\tabularnewline
			\cline{2-10} 
			& {\small{Truck}} & {\small{Sedan}} & {\small{HB}} & {\small{Truck}} & {\small{Sedan}} & {\small{HB}} & {\small{Truck}} & {\small{Sedan}} & {\small{HB}}\tabularnewline
			\hline 
			{\small{Weight (lbs)}} & {\small{6,050}} & {\small{3,216}} & {\small{2,496}} & {\small{6,800}} & {\small{3,190}} & {\small{2,295}} & {\small{9,900}} & {\small{2,097}} & {\small{3,620}}\tabularnewline
			\cline{1-1} 
			{\small{Length (inches)}} & {\small{207}} & {\small{195}} & {\small{162}} & {\small{229}} & {\small{189}} & {\small{154}} & {\small{233}} & {\small{179}} & {\small{160}}\tabularnewline
			\cline{1-1} 
			{\small{Height (inches)}} & {\small{70}} & {\small{58}} & {\small{60}} & {\small{76}} & {\small{58}} & {\small{59}} & {\small{77}} & {\small{58}} & {\small{58}}\tabularnewline
			\cline{1-1} 
			{\small{Width (inches)}} & {\small{78}} & {\small{73}} & {\small{67}} & {\small{80}} & {\small{72}} & {\small{67}} & {\small{80}} & {\small{72}} & {\small{68}}\tabularnewline
			\hline 
		\end{tabular}
		\par\end{centering}{\small \par}
	\caption{Vehicle dimensions (HB denotes Hatchback)}
	\label{fig:dim}
\end{table}

Every year, a new range of vehicles with different dimensions and weights are introduced. For instance, Ford manufactures trucks that weigh around 9,900 lbs. These kind of heavy vehicles might require two loading ramps for transportation. The vehicle dimensions for three types of vehicles from three different OEMs are shown in Table \ref{fig:dim}. The range of variations in the dimensions for the vehicle types indicates that it is imperative to find a good mix of vehicles for loading each auto-carrier so that the total number of auto-carriers required for deliveries is minimized. This necessitates a mathematical model to minimize the number of auto-carriers required for delivery of the vehicles. Other motivations include the increase in the number of vehicles being sold every year, fuel cost, and ever increasing variety of vehicles released in the market. 

The problem we address in the paper is described as follows: \textit{given a heterogeneous fleet of auto-carriers in a distribution center, a region consisting of a set of dealers, each requiring a set of vehicles, determine the route for the auto-carriers based on user's inputs and load the vehicles onto the auto-carriers to serve all the dealers at a minimum cost. Split deliveries are allowed, and the loads are not restricted by LIFO policy.}  We will refer to this problem as the auto-carrier transportation problem (ATP). This problem is known to be NP-hard \cite{Tadei2002}.
 
ATP consists of two underlying subproblems: routing the auto-carriers to deliver the vehicles to the dealers (routing subproblem) and finding a feasible load for an auto-carrier (loading subproblem). The solution to the loading subproblem depends on the sequence in which an auto-carrier makes its deliveries. In this paper, we generate a route for each auto-carrier and thereby fix the sequence of deliveries using a routing heuristic, \emph{i.e.}, we assume that the sequence of dealers that each auto-carrier delivers the vehicles is fixed a priori. Given this route sequence, we formulate the loading subproblem as a mixed integer linear program (MILP). We present a branch-and-price (B\&P) algorithm to optimally load the vehicles using the given heterogeneous fleet of auto-carriers. The algorithm generates a set of feasible loads. A feasible load includes the following information: the set of vehicles in the load, the auto-carrier onto which the vehicles are loaded, the vehicle-ramp assignment for the auto-carrier, position of the vehicle (whether it is loaded in the forward or reverse direction), slides and slide angle of each ramp in the auto-carrier. 

\subsection{Literature Review\label{sec:LitReview}}

ATP was first addressed in the literature by \cite{Agbegha1992} and \cite{Agbegha1998}. The authors describe the best practices followed by logistics companies in the U.S., and formulate the loading subproblem as a non-linear assignment problem. Subsequently, branch-and-bound algorithm is presented for the loading problem, however the routing subproblem is ignored. In the loading problem, an auto-carrier is modeled as a set of slots, and a loading network is introduced to impose a LIFO precedence among the slots. 
 
The authors of \cite{Tadei2002} present a case study of an Italian vehicle transportation company, and formulate ATP as a MILP. Planning horizon is multiple days, so the problem becomes complex. Due to this, they relax both the routing and loading problems. With regard to routing, the destinations are divided into multiple clusters, so the algorithm assigns the auto-carriers to the clusters. For loading, the vehicle lengths are approximated by equivalent constants and equated against the total length of an auto-carrier. Hence, the algorithm does not specify individual vehicle-ramp assignments. A greedy heuristic for the loading problem is developed in \cite{Miller2003}, and modeled an auto-carrier as vehicle with two flat loading platforms. The solution method assumes that the vehicles are always loaded straight on the platforms, and the loading is considered as a bin-packing problem with two bins.  

The work in \cite{Cuadrado2009} considers a real world auto-carrier distribution case in Venezuela, and develop a two-phase heuristic to determine a good sized fleet of auto-carriers based on a MILP formulation. Research in \cite{Lin2010} models the vehicle distribution in the U.S as a facility location problem, and presents a MILP formulation. However, the model does not explicitly consider loading and routing problems. Recently, the authors of \cite{Dell2014} propose an iterative local search algorithm for the routing, and enumerations techniques for the loading problem. LIFO policy is imposed to avoid reloading. We also refer the reader to \cite{iori2010routing} for a recent survey on loading and routing problems.  Stochastic version of auto-carrier loading problem is considered in \cite{saran}, and a special type of valid inequalities called Fenchel cutting planes. The details of Fenchel cutting planes can be found at \cite{beier2015stage}, \cite{venkatachalam2016sc}, and \cite{venkatachalam2016integer}.

We use B\&P approach for loading problem. To our knowledge, there is no work in the literature that develops an exact algorithm based on B\&P for the loading problem. A B\&P approach for solving an integer programming problem is similar to the conventional branch-and-cut approach but for the row generation procedure of branch-and-cut. In B\&P, column generation is used to solve the linear programs at each node of the branch-and-bound tree \cite{Barnhart1998}. This technique is especially useful when the number of decision variables in problem formulation increases exponentially with the size of the problem, and when explicit listing of all the columns becomes difficult. In such a scenario, columns are generated as needed at each node of the branch-and-bound search tree. The B\&P approach has been applied effectively to solve a variety of NP-hard problems such as assignment problem \cite{Savelsbergh1997}, shift scheduling \cite{Mehrotra2000}, vehicle routing problem (\cite{Dell2006}; \cite{Gutierrez2010}; \cite{Muter2014}), product line design \cite{Wang2009}, inventory routing problem \cite{Gronhaug2010}, service network design \cite{Andersen2011}, and joint tramp ship routing and bunkering \cite{Meng2015}. Readers interested in B\&P and column generation can refer to three excellent tutorials presented in \cite{Lubbecke2005}, \cite{wilhelm2001technical}, and \cite{Barnhart1998}.

A common feature among all the models presented in the literature for the loading problem is the usage of a variety of coefficients to model the dimensions of vehicles and auto-carriers. In practice, obtaining or estimating such coefficients from the past experience is not trivial as the OEMs continuously change the dimensions of the new vehicles every year. On the other side we have strict government regulations with regard to the total dimensions of a loaded auto-carrier, hence using exact parameters and calculations for the mathematical models become imperative from an optimization perspective, and also improves the practice. In summary, we propose a methodology to sequentially solve ATP; we also use exact dimensions of the vehicles and auto-carriers for modeling the loading problem. Using actual dimensions provide a means to replicate the methodology to suit other markets. Also, we do not impose LIFO restrictions on loading. 

\subsection{Objectives and contributions\label{sec:LitReview}}
In this paper, we describe a heuristic to generate a route based on user's inputs for the auto-carriers, and subsequently, we present an exact algorithm to optimally load auto-carriers for an auto logistics company. The contributions include a modeling framework for loading problem considering the actual dimensions of the vehicles, government and OEMs' restrictions, and loading sequence based on a fixed route sequence. Rather than restricting loading to LIFO scheme, the maximum number of reloads is a user input for the B\&P algorithm. Additionally, we provide the details to generate a route sequence. Furthermore, results from extensive computational experiments for loading problem using real-world instances are reported. The rest of the paper is structured as follows. Section \ref{sec:bp} provides the details of a heuristic to generate a route for the auto-carriers. We introduce notation and formulate the loading problem in Section \ref{sec:Formulation}, and describe a heuristic procedure for an initial feasible solution. Subsequently, we present B\&P algorithm for the loading problem. Section \ref{sec:compres} contains the computational study of the B\&P algorithm, and finally, Section \ref{sec:Conclusion} presents the concluding remarks and future research directions.

\section{Routing Heuristic\label{sec:bp}}
In this section, we present the details of the routing heuristic. The routing heuristic is used to determine the route, \emph{i.e.,} the sequence of deliveries to be made for each auto-carrier. We will refer to this sequence of deliveries as the routing sequence. The routing sequence is subsequently used by B\&P algorithm to generate loads for the auto-carriers.

\subsection{Routing heuristic\label{subsec:algo}}
As mentioned earlier, the ALC divides the dealers' locations into regions, and loads are built for each region. The user provides a source, a destination, an angle called the viewing angle, and an offset distance as inputs. The dealers' locations in the destination region are known a priori. Based on the user inputs, a polygon with four vertices is constructed. The details of the construction is elaborated in the Appendix. Then, the dealers' locations within the polygon are determined. The routing sequence for each auto-carrier is generated by solving a Hamiltonian path problem for the selected dealers' locations. The main steps of routing sequence are given in Figure \ref{BPIFSH}.

\begin{figure}[h]
	\begin{hanglist}
		\item \textbf{Routing sequence}
		\item \hrulefill
		
		{\it\bf \item Step 1. Initialization.} Obtain user inputs on source $dc$ and destination $d$ locations, viewing angle $\varphi$ and offset distance $\nu$. Let $la$ and $ln$ represent latitude and longitude for a location.
		
		{\it\bf \item Step 2. Calculate bearing angle.} Using the parameters $dc_{la}, dc_{ln}, d_{la},$ and $d_{ln}$, bearing angle $\theta$ is calculated.
		
		{\it\bf \item Step 3. Construct polygon.} Construct a polygon $\chi$ by computing three extreme points based on the values of $\theta$ and $\nu$.
		
		{\it\bf \item Step 4. Set of dealers' locations.} Form a set $S$ comprising of dealers' locations which lie inside the polygon $\chi$.		
		
		{\it\bf \item Step 5. Sequence generation.} Use Lin-Kernighan heuristic to solve a Hamiltonian path problem between $dc$ and $d$.
		\item \hrulefill
	\end{hanglist}
	\caption{Routing heuristic}
	\label{BPIFSH}
\end{figure}

In Step 1, the user chooses a source $dc$ and destination $d$ locations, viewing angle $\varphi$ and offset distance $\nu$. A source location $dc$ is ALC's distribution center, and a destination is a primary dealer's location where the load is to be built. Other parameters $\varphi$ and $\nu$ give the flexibility for the user to target the locations within a region for delivery. Using the latitude and longitude information for $dc$ and $d$, a bearing angle $\theta$ is calculated in Step 2. The bearing angle is defined as the angle made by the straight line between $dc$ and $d$ with respect to the geographical north. Based on the user's offset distance and bearing angle $\theta$, a polygon $\chi$ is constructed in Step 3. The set $S$ consists of all the dealer locations that lie inside the polygon $\chi$. The details of the formulae used from Step 2 to Step 4 are given in Appendix. Based on $dc$ and $d$, and the set of dealers' locations in the set $S$, a Hamiltonian path problem is solved using Lin-Kernighan heuristic \cite{LKH2000}. The output of Lin-Kernighan heuristic will provide the \textit{order} for each location $s \in S$. 

Another objective for the route sequence, that is commonly used as a performance indicator by the ALCs, is to improve the `percentage of perfect load' (PPL). This is explained as follows. Given the set of dealers' locations $S$, let $rm$ be the total miles traveled by an auto-carrier from the distribution center $dc$ to the last destination in the route sequence where the vehicles are delivered, $dt_s$ be the distance between a dealer location $s$ and ALC's distribution center, and $n_s$ be the number of vehicles delivered to the location $s \in S$. Then, `Optimum pay miles' is defined as $rm*|R(t)|$, where $R(t)$ is the set of ramps for the auto-carrier of type $t$, and `Tariff miles' is defined as $\sum_{s \in S}dt_s n_s$. PPL is defined as the ratio between `Tariff miles' and `Optimum pay miles.' ALC constantly looks forward to improve the performance indicator PPL. In a way it encourages the loads to a single dealer location or to multiple dealer locations which are closer to each other. Based on the selected dealers locations, Lin-Kernighan heuristic generates a route by minimizing the distance between them. In order to maximize PPL, the user will always select a dealer destination $d$ with high demand of vehicles or other locations with high demand which are very closer to $d$ so that `PPL' is maximized. This is a reason for our sequential approach to solve ATP. By careful selection of parameters for the construction of polygon, the optimal solution for ATP may not be far away from the one obtained by the presented methodology.

\section{Loading Algorithm\label{sec:Formulation}}
In this section we present a B\&P algorithm for the loading problem. B\&P algorithm constitutes of two mathematical models, a master problem (MP) and a pricing subproblem (LDP). LDP generates a feasible load, and MP is solved by adding the loads as columns at every node in the search tree. We first present the formulations for MP and LDP followed by B\&P algorithm. 

\subsection{Master problem} \label{sec:FormulationMaster} 
We first introduce notation to formulate MP. ALC has a fleet of $\bar{T}$ auto-carriers split into $T$ types and a set of vehicles $V$, indexed by $v$ to be delivered to dealers' locations. Identifying $t$ with an auto-carrier type, let $T^t_{\max}$ denote the maximum number of type $t$ auto-carriers available, and let $P(t)$ denote the set of feasible loads for an auto-carrier type $t \in T$, and $P(t)$ is indexed by $p(t)$. A feasible load includes the following information: the set of vehicles in the load and the vehicle-ramp assignment for auto-carrier type $t$. For each $p(t) \in P(t)$, we denote by $p_{{v}(t)}$, the set of vehicles in $p(t)$. Let $c^t$ denote the operating cost for an auto-carrier type $t$, and $c^{{v}}$ denote a penalty associated with a vehicle ${v} \in V$, if ${v}$ is not delivered to the assigned dealer. In practice, the penalty $c^{{v}}$ associated with the non-delivery of vehicle ${v}$ is very high. Let $x^t_{p(t)}$ be a binary variable, which indicates whether a feasible load $p(t)$ is used or not for an auto-carrier type $t$. Let $u^{{v}}$ be a binary variable, equal to $1$ if vehicle ${v}$ is not a part of any feasible load and $0$ otherwise. The MP is given as follows:

\begin{alignat}{2}
\text{MP1} = \min &\sum_{t \in T}\sum_{p(t) \in P(t)} c^t x^t_{p(t)} + \sum_{{v} \in {V}} c^{{v}} u^{{v}} \label{fs-objfun}
\end{alignat}

\hspace{2.5cm}subject to:\\
\begin{alignat}{2} \
\sum_{p(t) \in P(t)} x^t_{p(t)} \leq T^t_{\max} \qquad \forall \, t\in T, \label{capEqn-mp} \\
\sum_{t \in T} \sum_{p(t) \in P(t):{v} \subset p_{{v}(t)}} x^t_{p(t)} + u^{{v}} \geq  1  \qquad  \forall \, {v}\in {V}, \label{dmdEqn} \\
x^t_{p(t)} \in \{0,1\}, u^{{v}} \geq 0 \qquad \, \forall t \in T, \, p(t) \in P(t), {v}\in {V}. \label{vareqn}
\end{alignat}

\noindent Constraints \eqref{capEqn-mp} enforce the capacity limitation on the number of auto-carriers for each type $t$ that can be used for loading. Constraints \eqref{dmdEqn} ensure that every vehicle in the inventory is either loaded onto some auto-carrier for delivery or a penalty variable is triggered for its non-delivery to its dealer. In the objective function \eqref{fs-objfun}, we minimize the total operating costs for the auto-carriers and penalties for unsatisfied demand. The linear relaxation of MP, indicated hereafter by LP-MP, is same as MP with the binary restrictions on the $x^t_{p(t)}$ variables in \eqref{vareqn} relaxed, i.e., constraints \eqref{vareqn} are replaced as follows:

\begin{alignat}{2} \
x^t_{p(t)} \geq 0, \,\, u^{{v}} \geq 0 \qquad \, \forall t \in T, \, p(t) \in P(t), {v}\in {V}. \label{vareqn1}
\end{alignat}

In the following subsection, we formulate the pricing problem that generates feasible loads for MP.
%
%

\subsection{Pricing problem} \label{sec:FormulationPricing} 
\begin{figure}
	\centering
	\includegraphics[scale=0.6]{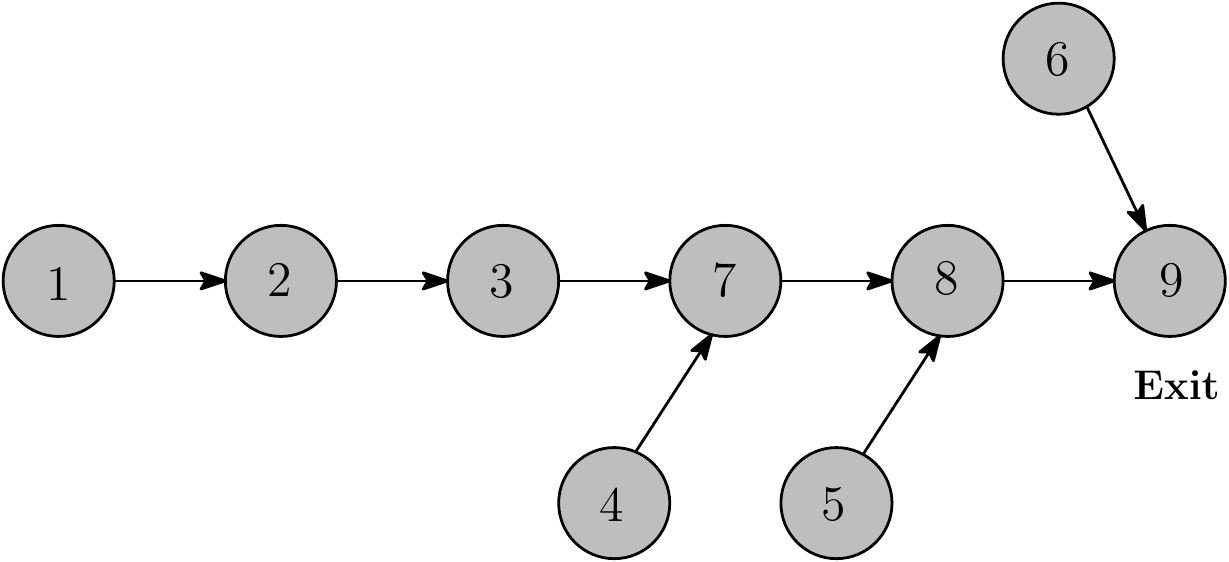}
	\caption{The loading/unloading graph for the auto-carrier in Figure \ref{fig:illustration}. Ramp $9$ is the exit ramp.}
	\label{fig:loadingsq}
\end{figure}

LDP is the pricing subproblem, and it finds a feasible load for the auto-carriers. As mentioned previously, a feasible load for an auto-carrier consists of a set of vehicles, and a precise vehicle-ramp assignment for each vehicle in an auto-carrier. The feasibility of a load is restricted by a number of constraints. The restrictions are due to the capacities and capabilities of the auto-carrier equipment, OEMs' preferences, and the government regulations that set limits on the length, height and weight of a loaded auto-carrier. We develop each of these constraints after introducing our notation. As shown in Figure \ref{fig:illustration}, an auto-carrier of type $t \in T$ consists of a set of ramps $R(t)$ indexed by $i$. For a given auto-carrier of type $t$ and its corresponding $R(t)$, we define the following subsets of $R(t)$:  
\begin{alignat}{2}\notag
R_{\mathcal{U}}(t) \equiv&\text{ set of ramps in the upper deck of an auto-carrier of type $t$ (e.g., ramps 1, 2, 4, 5 and 6 in}\notag \\
&\text{Figure \ref{fig:illustration}), and} \notag \\
R_{\mathcal{L}}(t) \equiv&\text{ set of ramps in the lower deck of an auto-carrier of type $t$ (e.g., ramps 3, 7, 8 and 9 in }\notag\\
&\text{Figure \ref{fig:illustration}).} \notag
\end{alignat}
For each ramp $i\in R_{\mathcal{U}}(t)$, we also define:
\begin{alignat}{2}\notag
r_L(i) \equiv& \text{ ramp at the lower deck for a given ramp $i \in  R_{\mathcal{U}}(t)$.} \notag
\end{alignat}
For a given auto-carrier of type $t$ and its corresponding ramp set $R(t)$, we define the following collection of subsets of $R(t)$ (the superscript $c$ in the following notation indicate that we are referring to a collection of sets):
\begin{alignat}{2}\notag
R^c_{\mathcal{H}}(t) \equiv& \text{ collection of sets of ramps that limit the height of a loaded auto-carrier of type $t$.} \\ \notag
& \text{ (e.g., $R_{\mathcal{H}}^c(t) = \{\{1\}, \{2,3\}, \{4,7\}, \{5,8\}, \{6,9\}\}$ for auto-carrier in Figure \ref{fig:illustration})}. \\ \notag
R^c_{\mathcal{L}}(t) \equiv& \text{ collection of sets of ramps in an auto-carrier of type $t$, that have a limitation on the}\\ \notag
& \text{allowed length (e.g., $R^c_{\mathcal{L}}(t) = \{\{1,2\}, \{3\}, \{4,5,6\}, \{7,8,9\}\}$ for auto-carrier in Figure \ref{fig:illustration})}. \\ \notag	
R^c_{\mathcal{SP}}(t) \equiv&\text{ collection of sets of ramps where each set of ramps can be combined to accommodate a}\\ \notag
& \text{ single vehicle in an auto-carrier of type $t$. (e.g., $R_{\mathcal{SP}}^c(t) = \{\{4,5\}, \{7,8\}, \{5,6\}\}$ for} \\ \notag
& \text{ an auto-carrier in Figure \ref{fig:illustration}).} \notag	
\end{alignat}
We refer to $R^c_{\mathcal{SP}}(t)$ as the collection of \emph{split ramps} for an auto-carrier of type $t$. 
Another aspect of LDP is a loading/unloading sequence for each ramp in an auto-carrier. The sequence in which the vehicles are loaded or unloaded from the ramps of an auto-carrier is unique and is dictated by a loading/unloading graph (see Figure \ref{fig:loadingsq}). The vertices of the graph represent the ramps of an auto-carrier, and they denote the set of ramps which should be empty or unloaded before a vehicle in a particular ramp is unloaded. To illustrate, if a vehicle in ramp 4 should be unloaded then the vehicles in the ramps 7, 8, and 9 should be unloaded or empty. A ramp in the graph without any outgoing edges is referred as an exit ramp (9 in Figure \ref{fig:loadingsq}). For every ramp $i \in R(t)$ of an auto-carrier of type $t$, the path from $i$ to the exit ramp is unique, and this gives an unloading sequence for the vehicle in the ramp $i$. Hence, for every ramp $i\in R(t)$, we define $q^t(i)$ as the ordered set of ramps in which the vehicles should be unloaded before unloading the vehicle in ramp $i$.
\begin{alignat}{2}\notag
q^t(i) \equiv&\text{ ordered set of ramps providing the unloading sequence of ramp $i$ in an auto-carrier of type $t$.} \notag
\end{alignat}

For each auto-carrier, we also assume a fixed route, i.e, an ordered sequence of dealer locations that each auto-carrier is going to visit. The set $S$ is an \textit{ordered} set of dealers' locations is obtained by the route heuristics explained in the previous section. For any $s_1,s_2 \in S$, $s_1 \preceq s_2$ represents that location $s_1$ will be visited by the auto-carrier before visiting the location $s_2$, and $\succeq$ represents vice-versa.  We now define for each $v \in V$, $lo(v)$ as the dealer location. 
\begin{alignat}{2}\notag
lo(v) \equiv&\text{  dealer location for a vehicle $v \in V$}, lo(v) \in S. \notag
\end{alignat}
Now, we will describe various configurations in which a vehicle can be loaded on a ramp in an auto-carrier. Suppose that a vehicle $v$ is loaded on a ramp $i \in R(t) \cup R^c_{\mathcal{SP}}(t)$ then $v$ can either be positioned in the same direction as that of an auto-carrier (vehicles positioned in the ramps 2, 3, 4, 5, 6 and 9 in Figure \ref{fig:illustration}) or in the opposite direction (vehicles positioned in the ramps 1, 7 and 8 in Figure \ref{fig:illustration}). Once $v$ is positioned on a ramp $i$, the ramp can either slide in the forward direction (e.g., ramps 5 and 8 in Figure \ref{fig:illustration}) or in the reverse direction (e.g., ramps 2, 4, 6 and 9 in Figure \ref{fig:illustration}). The amount of forward or backward slide is restricted by a maximum allowable angle of slide for each ramp. We now define the following sets:
\begin{flalign}
&J \equiv\text{ set of possible positions for a vehicle in a ramp, indexed by $j$,}\notag \\
&L \equiv\text{ set of possible slides for a ramp in an auto-carrier, indexed by $\ell$,}\notag \\
&M \equiv\text{ discrete set of allowable slide angles for a ramp in an auto-carrier, indexed by $m$, and} \notag \\
&\boldsymbol{\alpha}^t(R,J,V,L,M)\equiv\text{ set of all possible configurations $(i,j,v,l,m)$, where $i \in R$, $j\in J$, $v\in V$,}\notag \\
&\qquad\qquad\qquad\qquad\text{$\ell \in L$, and $m \in M$ in an auto-carrier of type $t$.}& \notag
\end{flalign}

\begin{figure}
	\centering
	\includegraphics[scale=0.55]{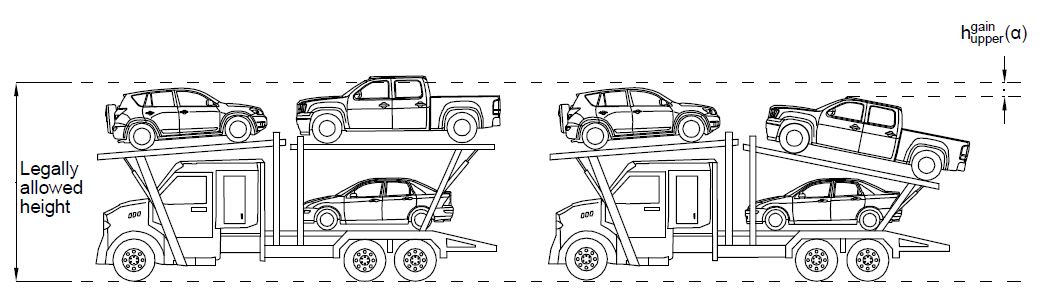}
	\caption{$h^{\text{gain}}_{\text{upper}}(\alpha)$ for a particular configuration $\alpha$. The height of the loaded auto-carrier exceeds the legally allowed height when the ramp does not slide. When the ramp slides at an angle, a gain in height is obtained and the height restriction is met.}
	\label{fig:htadv}
\end{figure}

It should be noted that the set $\boldsymbol{\alpha}^t(R,J,V,L,M)$ is constructed based on the inclusion and exclusions preferences for vehicle-ramp assignment from OEMs, and also, the set excludes the vehicles from the ramps due to incompatible dimensions, i.e., if a wheel base of a vehicle is larger than a ramp base. We now define a set of parameters whose values are either obtained from the auto-carrier's or vehicle's specifications or restrictions imposed by the governmental agencies for a loaded auto-carrier. 
\begin{alignat}{2}\notag
l^t_{\max}(\mathcal{L}) \equiv&\text{ maximum allowable loaded length for $\mathcal{L} \in R_{\mathcal{L}}^c(t)$ for an auto-carrier of type $t$, } \\ \notag
h_{v}				\equiv&\text{ height of the vehicle $v \in V$,} \\ \notag
w_{v}				\equiv&\text{ weight of the vehicle $v \in V$,} \\ \notag
h_{\max}			\equiv&\text{ maximum allowable height (loaded height) for an auto-carrier,} \\ \notag
w^{\text{steer}}_{\max}	\equiv&\text{ maximum allowable weight at the steering axle for an auto-carrier,}\\ \notag
w^{\text{drive}}_{\max}	\equiv&\text{ maximum allowable weight at the drive axle for an auto-carrier,}\\ \notag
w^{\text{tandem}}_{\max}	\equiv&\text{ maximum allowable weight at the tandem axle for an auto-carrier,}\\ \notag
w^{\text{total}}_{\max}\equiv&\text{ maximum allowable total load weight for an auto-carrier.} \notag
\end{alignat}
\begin{figure}
	\centering
	\includegraphics[scale=0.35]{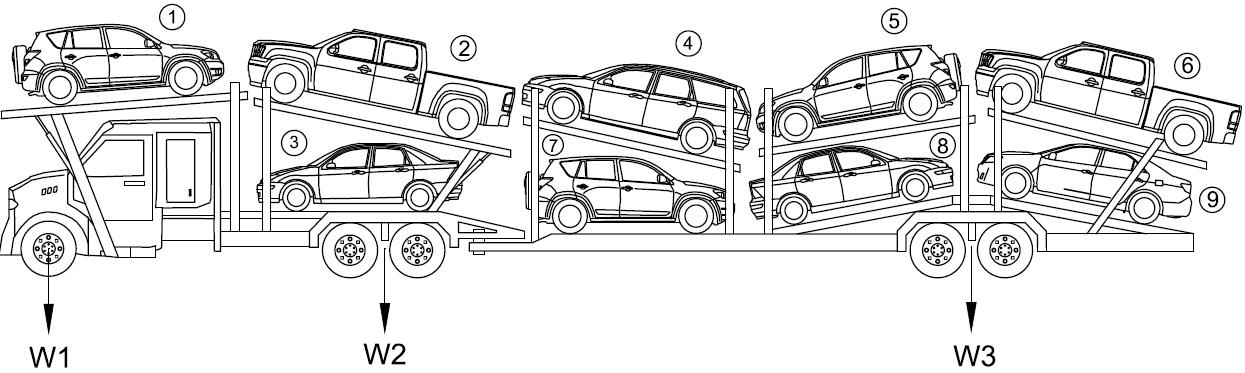}
	\caption{$W_1$, $W_2$, and $W_3$ are the steer, drive, and tandem axles, respectively.}
	\label{fig:axles}
\end{figure}
Figure \ref{fig:axles} shows the steer, drive, and tandem axles for an auto-carrier. We now introduce parameters that are derived using vehicle and auto-carrier specifications. The reader is referred to the Appendix for the details of the calculations for each of these parameters. 
\begin{alignat}{2}\notag
l^{v}_m 			\equiv&\text{ effective length of a vehicle $v \in V$ when it is loaded on a ramp inclined at an angle $m$,}& \\ \notag
h^{\text{gain}}_{\text{upper}}(\alpha) \equiv&\text{ height gain for a configuration $\alpha \in \boldsymbol{\alpha}^t(R_{\mathcal{U}}(t),J,V,L,M)$,} \\ \notag
h^{\text{max}}_{\text{lower}}(\alpha) \equiv&\text{ maximum allowable slide for a configuration $\alpha \in \boldsymbol{\alpha}^t(R_{\mathcal{L}}(t),J,V,L,M)$,} \\ \notag
w^{\text{steer}}(\alpha) \equiv&\text{ weight contributed to steer axle by configuration $\alpha \in \boldsymbol{\alpha}^t(R(t)\cup R^c_{\mathcal{SP}}(t),J,V,L,M)$,} \\ \notag
w^{\text{drive}}(\alpha) \equiv&\text{ weight contributed to drive axle by configuration $\alpha \in \boldsymbol{\alpha}^t(R(t)\cup R^c_{\mathcal{SP}}(t),J,V,L,M)$,} \\ \notag
w^{\text{tandem}}(\alpha) \equiv&\text{ weight contributed to tandem axle by configuration $\alpha \in \boldsymbol{\alpha}^t(R(t)\cup R^c_{\mathcal{SP}}(t),J,V,L,M)$.} \\ \notag
\end{alignat}
The parameters $h^{\text{gain}}_{\text{upper}}(\alpha)$ and $h^{\text{max}}_{\text{lower}}(\alpha)$ are derived based on the dimensions of the vehicles and ramps of an auto-carrier. The parameter $h^{\text{gain}}_{\text{upper}}(\alpha)$ is the gain in height due to the sliding of upper ramp by an angle, and $h^{\text{max}}_{\text{lower}}(\alpha)$ indicates the maximum possible slide the auto-carrier can accommodate for a vehicle in the upper deck due to the contour of the vehicle in the lower deck. Figure \ref{fig:htadv} shows an illustration of $h^{\text{gain}}_{\text{upper}}(\alpha)$ for a particular configuration $\alpha$.
A loaded auto-carrier visits each location and unloads the vehicles that are to be delivered to that dealer. While unloading a vehicle $v_1$ from an auto-carrier, a violation is said to occur if a vehicle to be delivered to another dealer location has to be unloaded from the auto-carrier in order to unload $v_1$. Violations can be avoided by using a LIFO policy while loading, however, this may increase the number of auto-carriers required for delivery. In practice, ALC imposes an upper bound on the number of violations that an auto-carrier can incur. Hence, we define
\begin{alignat}{2}\notag
\text{v}_{\max} \equiv&\text{ maximum number of unload violations that an auto-carrier can incur. }  \notag
\end{alignat}

\noindent{\bf Decision variables}\\
\\
We now define the decision variables used to formulate LDP. For every auto-carrier of type $t$ and $\alpha \in \boldsymbol{\alpha}^t(R(t)\cup R^c_{\mathcal{SP}}(t),J,V,L,M)$, let ${y}(\alpha)$ be a binary variable, equal to $1$ if the configuration $\alpha$ is used and $0$ otherwise. 
\begin{align}
{y}(\alpha) &=\left\{ \begin {array}{clcr}
1, & \mbox{if configuration $\alpha\in \boldsymbol{\alpha}^t(R(t)\cup R^c_{\mathcal{SP}}(t),J,V,L,M)$ is in use,} \\  
0, & \mbox{otherwise.}
\end {array}
\right. \notag \\ \notag
\end{align}
For every auto-carrier of type $t$ and $i\in R^c_\mathcal{SP}(t)$, let ${sp}_i$ be a binary variable defined as follows:

\begin{align}
{sp}_{i} &=\left\{ \begin {array}{clcr}
1, & \mbox{if the split ramp $i\in R^c_\mathcal{SP}(t)$ is used}, \\
0, & \mbox{otherwise.}
\end {array}
\right. \notag \\ \notag
\end{align}
We also define a binary variable ${z}_{is}$ to denote the ramp and location relationship, and an integer variable $u_{is}$ to indicate the number of violations for the vehicle in ramp $i$ based on the location $s$. 
\begin{align}
{z}_{is} &=\left\{ \begin {array}{clcr}
1, & \mbox{if any vehicle $v\in V$ with $lo(v)=s$ is assigned to ramp $i\in R(t)\cup R^c_{\mathcal{SP}}(t)$}, \\
0, & \mbox{otherwise.}
\end {array}
\right. \notag \\
u_{is} &= \text{ number of violations for the vehicle in ramp $i$, delivered to location $s$.}\notag 
\end{align}

We now define a continuous variable:
\begin{align}
{h}_{i \ell } \, \, \equiv &\text{ adjusted height gain in the ramp $i \in R(t)$ at slide position $\ell$.} \notag 
\end{align}

Finally, we define a dual value ${d}^v$ for each constraint in \eqref{dmdEqn} of MP, and $r^v$ is the potential revenue for delivering the vehicle $v$ to its dealer location. The revenue depends on the vehicle type and the location of the dealer. 

\subsection{Model Formulation}
A MILP formulation for LDP is presented in this section. LDP for each auto-carrier of type $t$ (LDP($t$)) is presented as follows.

\begin{alignat}{2}
\text{LDP1}(t) = \min (1 - \sum_{\alpha \in \boldsymbol{\alpha}^t(R(t)\cup R^c_{\mathcal{SP}}(t),J,V,L,M)} ({d}^v+r^v) {y}(\alpha)) \label{ss-objfun}
\end{alignat}

Each configuration of $\alpha$ is a unique combination $(i,j,v,\ell,m)$ for an auto-carrier of type $t$, where $i$ is a ramp (single or split); $v$ is the vehicle on ramp $i$; $j$ is the position of vehicle $v$ on $i$; $\ell$ is the slide of ramp $i$; and $m$ is the slide angle of ramp $i$. \\

\noindent{\bf Constraints}\\
\\
The objective function is optimized over the set of feasible solutions described by the following sets of constraints: \\

\noindent{\bf Height constraints}
\begin{figure}
	\centering
	\includegraphics[scale=0.75]{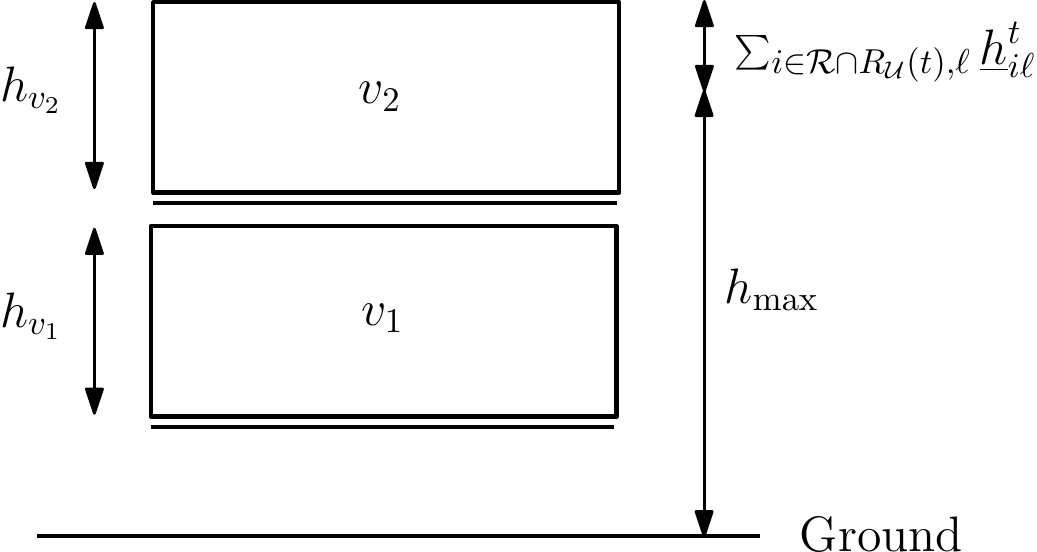}
	\caption{Height constraints in equation \eqref{ht_c1}.}
	\label{fig:htconst}
\end{figure}
\begin{alignat}{2}
\sum_{
	\substack{\alpha \in  \boldsymbol{\alpha}^t(R,J,V,L,M)}}  
h_v {y}(\alpha) - 
\sum_{i \in R \cap R_{\mathcal{U}}(t), \ell} {h}_{i \ell }
\leq h_{\max} \qquad  &\forall \, R \in R^c_{\mathcal{H}}(t),  \label{ht_c1} \\ 
{h}_{i \ell } \leq 
\sum_{\substack{\alpha \in  \boldsymbol{\alpha}^t(\{i\},J,V,\{\ell\},M)}} 
h^{\text{gain}}_{\text{upper}}(\alpha) {y}(\alpha)
\qquad  &\forall \, i \in R_{\mathcal{U}}(t), \,  \ell \in L, \label{ht_c2} \\
{h}_{i \ell } \leq 
{h}_{r_L(i)  \ell }
\qquad  &\forall \, i \in R_{\mathcal{U}}(t), \,  \ell \in L,  \label{ht_c3} \\
{h}_{i  \ell } \leq 
\sum_{\substack{\alpha \in  \boldsymbol{\alpha}^t(\{i\},J,V,\{\ell\},M)}} 
h^{\text{max}}_{\text{lower}}(\alpha) {y}(\alpha)
\qquad  &\forall \, i \in R_{\mathcal{L}}(t), \,  \ell \in L. \label{ht_c4} 
\end{alignat}

\noindent Constraints \eqref{ht_c1} enforce the maximum height limitation for a pair of vehicles in the lower and upper deck for an auto-carrier of type $t$ (see Figure \ref{fig:htconst}). The term $\sum_{i \in R\cap R_{\mathcal{U}}(t), \ell} {h}^{t}_{i \ell }$ takes value based on the pair of vehicles considered, and represents the total height that has to be reduced by sliding the ramps. Also, this height is restricted by the contour of vehicles $v_1$ and $v_2$ (constraints \eqref{ht_c2}--\eqref{ht_c4}). Each element of the set $R \in R^c_{\mathcal{H}}(t)$ represents a set of ramps, one in a lower and another in the corresponding upper deck, which are vertically aligned for an auto-carrier of type $t$, and constitutes the maximum height limitation. The variable ${h}_{i \ell }$ is the height gain achieved by sliding ramp $i$ at an angle $\ell$. For ramp $i$ in the upper deck, the constraints in \eqref{ht_c2} and \eqref{ht_c3} upper bound the variables ${h}_{i \ell }$. The upper bound in constraint \eqref{ht_c2} is due to the vehicle in the ramp $i$, and the upper bound in \eqref{ht_c3} is due to the vehicle in the ramp below $i$, \emph{i.e.}, $r_L(i)$. Similarly, constraints \eqref{ht_c4} bound the variables ${h}_{i\ell }$ based on the vehicle in ramp $i$ in the lower deck. \\

\noindent{\bf Weight constraints}
\begin{alignat}{2}
\sum_{\substack{\alpha \in \boldsymbol{\alpha}^t(R(t)\cup R^c_{\mathcal{SP}}(t),J,V,L,M)}}
w^{\text{steer}}(\alpha) {y}(\alpha) 
&\leq w^{\text{steer}}_{\max}, \label{wt_c1} \\   	
\sum_{\substack{\alpha \in \boldsymbol{\alpha}^t(R(t)\cup R^c_{\mathcal{SP}}(t),J,V,L,M)}}
w^{\text{drive}}(\alpha) {y}(\alpha) 
&\leq w^{\text{drive}}_{\max}, \label{wt_c2} \\   	
\sum_{\substack{\alpha \in \boldsymbol{\alpha}^t(R(t)\cup R^c_{\mathcal{SP}}(t),J,V,L,M)}}
w^{\text{tandem}}(\alpha) {y}(\alpha) 
&\leq w^{\text{tandem}}_{\max}, \label{wt_c3} \\   	
\sum_{\substack{\alpha \in \boldsymbol{\alpha}^t(R(t)\cup R^c_{\mathcal{SP}}(t),J,V,L,M)}}
w_v {y}(\alpha) 
&\leq w^{\text{total}}_{\max}. \label{wt_c4}  
\end{alignat}
%
\noindent The total weight for each axle of an auto-carrier of type $t$ is restricted in the weight constraints. The three axle weights, namely steer axle weight, drive axle weight, and tandem axle weight are restricted in the constraints \eqref{wt_c1}, \eqref{wt_c2}, and \eqref{wt_c3}, respectively. Constraint \eqref{wt_c4} restricts the total weight for an auto-carrier. 

\vspace{0.5cm}
\noindent{\bf Length constraints}
\\
\begin{alignat}{1}
\sum_{v,m}\left( l^v_m  
\sum_{\substack{\alpha \in \boldsymbol{\alpha}^t(\mathcal{L},J,\{v\},L,\{m\})}} {y}(\alpha)
\right)\leq l_{\max}^t(\mathcal{L}) 
\qquad  \forall \, \mathcal{L} \in R^c_{\mathcal{L}}(t).  \label{ln_c1} 
\end{alignat}

\noindent Constraints \eqref{ln_c1} restrict the total length for each ramp set $\mathcal{L} \in R^c_{\mathcal{L}}(t)$ for an auto-carrier of type $t$.\\

\noindent{\bf Split Ramp constraints}
\begin{alignat}{1}
\sum_{\substack{\alpha \in \boldsymbol{\alpha}^t(\{i\},J,V,L,M)}} {y}(\alpha) \leq {sp}_i  
\qquad  &\forall \, i \in R^c_{\mathcal{SP}}(t),  \label{sp_c1} \\
\sum_{\alpha \in \boldsymbol{\alpha}^t(i\cap R(t),J,V,L,M)} {y}(\alpha) \leq |i|(1-{sp}_i)
\qquad  &\forall \, i \in R^c_{\mathcal{SP}}(t).  \label{sp_c2} 
\end{alignat}

\noindent Constraints \eqref{sp_c1} and \eqref{sp_c2} ensure that when a vehicle is loaded on a split ramp by combining two or more ramps, then assignments of vehicles are not made to the individual ramps that are combined to form the split ramp and vice versa. 

\vspace{0.5cm}
\noindent{\bf Assignment constraints}
\begin{alignat}{1}
\sum_{\substack{\alpha \in \boldsymbol{\alpha}^t(\{i\},J,V,L,M)}} {y}(\alpha) \leq 1
\qquad  &\forall \, i\in R(t) \cup R^c_{\mathcal{SP}}(t),  \label{as_c1} \\
\sum_{\substack{\alpha \in \boldsymbol{\alpha}^t(R(t) \cup R^c_{\mathcal{SP}}(t),J,\{v\},L,M)}} {y}(\alpha) \leq 1
\qquad  &\forall \, v\in V.  \label{as_c2} 
\end{alignat}

\noindent Constraints \eqref{as_c1} and \eqref{as_c2} ensure that every vehicle is assigned to only one ramp and each ramp has only one vehicle. 

\vspace{0.5cm}
\noindent{\bf Violation constraints}
\begin{alignat}{1}
{z}_{is} \geq 
\sum_{\substack{\alpha \in \boldsymbol{\alpha}^t(\{i\},J,\{v\in V:lo(v)=s\},L,M)}}	{y}(\alpha) 
\qquad  &\forall \, i\in R(t) \cup R^c_{\mathcal{SP}}(t), \, \, s\in S, \label{vi_c1} \\
u_{is} \geq \sum_{\substack{\alpha \in \boldsymbol{\alpha}^t(q^t(i),J,\{v\in V:lo(v) \succeq s\},L,M)}} {y}(\alpha) 
- |R(t)|(1-{z}_{is}) 
\qquad  &\forall \, i\in R(t) \cup R^c_{\mathcal{SP}}(t), \, \, s\in S, \label{vi_c2} \\
\sum_{i \in R(t) \cup R^c_{\mathcal{SP}}(t), s \in S} u_{is} \leq \text{v}^t_{\max}. \qquad \label{vi_c3} 	
\end{alignat}

\noindent Constraints \eqref{vi_c1} assign value for the variables ${z}_{is}$ based on vehicle $v$'s location in the ramp $i$. Constraints \eqref{vi_c2} account the violations for the ramp $i$ based on the route sequence of the assigned vehicle $v$ and the loading sequence $q^t(i)$ for the ramp $i$. Subsequently, constraint \eqref{vi_c3} limits the number of violations that can occur.\\

Height, weight, and length constraints ensure that the generated load satisfies the regulations imposed by the governmental authorities. Assignment constraints avoid duplicate allotments for a ramp or a vehicle, and finally, violation constraints restrict the number of reloads for an auto-carrier.\\

To corroborate the performance of the B\&P algorithm, an equivalent aggregated formulation (EP) of the loading problem is developed and solved using a direct solver. The formulation for EP is briefly explained in this section. Let $K(t)$ be the set of auto-carriers for $t \in T$, and variable ${y}(\alpha)$ is added with indices $i$ and $t$ such that the variable ${y}^t_i(\alpha)$ indicates $i^{th}$ auto-carrier of type $t$. Let ${p}^t_i$ denote a binary variable indicating whether the auto-carrier ${y}^t_i(\alpha)$ is used or not. The EP model is given as follows.

\begin{alignat}{2}
\text{EP1} = \min \sum_{ t\in T, \, \, i \in K(t) } c^t {p}^t_i  \label{ep-objfun}
\end{alignat}
\hspace{2.5cm}subject to:\\
\begin{alignat}{2} \
\sum_{\substack{\alpha \in \boldsymbol{\alpha}^t(R(t) \cup R^c_{\mathcal{SP}}(t),J,V,L,M)}} {y}_i^t(\alpha) \leq |R(t)| {p}^t_i
\qquad  \forall \, t\in T, \, i \in K(t), \label{ep-eqn1} 
\end{alignat}

\begin{alignat}{2} \
\sum_{i \in K(t)} {p}^t_i \leq T^t_{\max} \qquad \forall \, t\in T, \label{ep-eqn2}
\end{alignat}

\begin{alignat}{2} \
\sum_{\substack{i \in K(t) ,\alpha \in \boldsymbol{\alpha}^t(R(t) \cup R^c_{\mathcal{SP}}(t),J,v,L,M)}} {y}_i^t(\alpha) \geq 1
\qquad  &\forall \, v\in V.  \label{ep-eqn3} 
\end{alignat}

We minimize the total cost required for delivering the vehicles. In constraints \eqref{ep-eqn1}, variable ${p}^t_i$ indicates whether an $i^{th}$ auto-carrier of type $t$ is used or not. The number of auto-carriers for each type $t$ is capacitated by constraint \eqref{ep-eqn2}. Constraints \eqref{ep-eqn3} assure that every vehicle $v$ is loaded in at least one auto-carrier $i$ of type $t$. Constraints \eqref{ht_c1} to \eqref{vi_c3} are included in EP model with indices $i$ and $t$ for the variable ${y}(\alpha)$.

\subsection{Branch-and-price algorithm\label{subsec:algo}}
MP in Section \ref{sec:FormulationMaster} has a very large number of decision variables so we develop a B\&P approach to solve the problem. For an optimal solution to LP-MP, we need to explicitly consider the entire set of possible variables by $P(t)$ which is practically very exhaustive. Hence, we construct another model called a restricted master problem (RMP) which relaxes the LP-MP using a smaller set $\bar{P}(t)$. Columns will be added to $\bar{P}(t)$ by solving a pricing problem, which in our case will be LDP1$(t)$ for an auto-carrier of type $t$. LDP1$(t)$ will be solved using the dual values obtained from RMP, and the generated load will be added to the set $\bar{P}(t)$ if it can improve the objective solution of the RMP further. The B\&P approach guarantees an optimal solution by generating columns at each node within the branch-and-bound search tree. In the following subsection, we discuss the main steps involved in the B\&P algorithm.

\begin{figure}
	\begin{hanglist}
		\item \textbf{Initial feasible solution heuristic}
		\item \hrulefill
		
		{\it\bf \item Step 1. Initialization.}  Initialize $d^v$ to a positive number.
		
		{\it\bf \item Step 2. Solve pricing problem.} Arbitrarily choose a $t$ such that $\left\vert \bar{P}(t)\right\vert < T^t_{max}$ and solve the pricing problem LDP1($t$), and let the solution be $p(t)$.
		
		{\it\bf \item Step 3. Update columns.} Add the solution to $\bar{P}(t)$, $\bar{P}(t)$ = $\bar{P}(t) \cup p(t)$. Set $V = V - v, \forall v \subset p(t)$.
		
		{\it\bf \item Step 4. Stop criterion. } If $V =\phi$ then Stop, else go to Step 2.		
		
		\item \hrulefill
	\end{hanglist}
	\caption{Initial feasible solution heuristic}
	\label{BPIFSH}
\end{figure}

\subsubsection{Algorithm\label{subsec:algo}}

An initial solution to the problem is constructed using initial feasible solution (IFS) heuristic. IFS is used to generate an initial set of columns for $\bar{P}(t)$ then the B\&P algorithm is used. The heuristic is given in Figure \ref{BPIFSH}.


In Step 1 we initialize the dual parameter $d^v$ to a positive number. We solve the pricing problem LDP1($t$) in Step 2 for an arbitrary auto-carrier of type $t$ such that there is enough capacity available for loading. Based on the output from the pricing problem LDP1($t$), the columns are added to the respective $\bar{P}(t)$. In Step 4 we check whether all the vehicles are loaded, if not then algorithm is directed to Step 2, otherwise, terminate the algorithm. By using IFS, we have set of initial loads for all the vehicles in the set $V$.

\begin{figure}
	\begin{hanglist}
		\item \textbf{Branch-and-price algorithm}
		\item \hrulefill
		
		{\it\bf \item Step 1. Initialization.}  Run the routing heuristic and generate route sequence $lo(v), \forall v \in V$. Set upper bound UB=$+\infty$, lower bound LB=$-\infty$, and let LDP1$(t) = 0 , \forall t \in T$, $\epsilon >0 $ as tolerance, and let $\hat{x}^t_{p(t)}$ be the solution for RMP. A node $n$ represents an instance of RMP, RMP-Obj($n$) represents objective value for RMP for the node $n$, $\Gamma$ represents a list with nodes, and a root node $n_0$ is formed based on the columns generated from IFS, and $n_0$ is added to $\Gamma$. Set node $n_0$ as $n$.
		
		{\it\bf \item Step 2. Solve pricing problem.} For every $t \in T$, solve the respective LDP($t$) based on the dual values of node $n$ as a MILP, and update LDP1$(t)$.
		
		{\it\bf \item Step 3. Pricing problem selection.} If LDP1$(t) \geq 0, \forall t \in T$, then goto Step 6. Otherwise, select a pricing problem LDP1$(t')$ such that $t' = \argmin_{t} \{\text{LDP1}(t)\}$, and add the columns based on the solution of LDP1($t'$) to $\bar{P}(t')$ for the node $n$, $\bar{P}(t') = \bar{P}(t') \cup p(t')$. This choice of subproblem is called `best-first' policy.
		
		{\it\bf \item Step 4. Solve master problem and branching.} Solve RMP for node $n$ and update the dual values. If there is any $\hat{x}^t_{p(t)} \notin \mathbb{Z} \, \, \, \, \forall \, \, t \in T, \, \,$  and LB < RMP-Obj($n$), then create two nodes $n_1$ and $n_2$ and copy all the information for the nodes from the node $n$. Choose a variable $x^{t^*}_{p(t)^*}$ such that $\argmax_{t^*,p(t)^*} \{\lvert \hat{x}^{t^*}_{p(t)^*} -1 \rvert : \hat{x}^{t^*}_{p(t)^*} \notin \mathbb{Z}\} $.  Add the following constraints $x^{t^*}_{p(t)^*} \leq \lfloor\hat{x}^{t^*}_{p(t)^*} \rfloor$ and $x^{t^*}_{p(t)^*} \geq \lceil \hat{x}^{t^*}_{p(t)^*} \rceil $ to the nodes $n_1$ and $n_2$, respectively. Add the nodes to the list $\Gamma$. Set LB = RMP-Obj($n$). If $\hat{x}^t_{p(t)} \in \mathbb{Z} \, \, \, \, \forall \, \, t \in T, \,\, p(t) \in \,\, \bar{P}(t) $ and UB > RMP-Obj($n$), then store the solution $\hat{x}^t_{p(t)}$ and set UB = RMP-Obj($n$).		
		
		{\it\bf  \item Step 5. Node fathoming.}  Remove any node $n^*$ from the list $\Gamma$ such that RMP-Obj($n^*$) $<$ LB or RMP-Obj($n^*$) $>$ UB. 		
		
		{\it\bf  \item Step 6. Stop criterion.} If $\Gamma = \phi $ or UB-LB $\leq \epsilon$, then Stop. Otherwise, use depth-first criterion to choose a node and set the node as $n$. Go to Step 2.		
		
		\item \hrulefill
	\end{hanglist}
	\caption{Branch-and-price algorithm}
	\label{BP}
\end{figure}

We now outline the main components of B\&P algorithm to compute optimal solutions for the loading problem. In Step 1 we use the routing and IFS heuristics to generate a route sequence and initial set of feasible loads, respectively.  An instance of RMP is a node $n$ and a list $\Gamma$ is used to store nodes. A root node is constructed using the outputs from IFS and added to the list $\Gamma$. Based on the dual values $d^v$ of node $n$, the objective function of LDP1($t$) is updated and solved as a MILP. If the objective functions for all the auto-carrier types are non-negative then stop criterion at Step 6 are checked. Otherwise, choose the pricing problem with the least objective value so that this could provide a maximum improvement for the objective function of RMP, and add the solution vector as a column to the node $n$. In Step 4, RMP is solved for node $n$, and if there exists a solution for a variable which is non-integral then two nodes are created using current node's information and a bisection of solution space is done based on the variable with non-integral value and added to the list. If the solution is integral for all the variables, then upper bound is checked and updated accordingly. Based on the updated lower and upper bounds, the nodes in the list $\Gamma$ are fathomed in Step 5. In Step 6, we evaluate the stopping criterion, and if the list $\Gamma$ is empty or the difference between lower and upper bounds is within a given tolerance then the algorithm is terminated. Otherwise, the algorithm continues from Step 2 with the updated information.

\section{Computational Results\label{sec:compres}}
In this section, we present the computational results of B\&P algorithm for the loading problem. Solution method for routing problem is predominantly based on well known Lin-Kernighan heuristic. Hence, we omit the computational study for the routing problem. The objectives of the computational study for the B\&P algorithm are \textit{(i)} test the efficacy of the procedure using real world instances \textit{(ii)} benchmark the performance against a holistic model, and \textit{(iii)} evaluate the solution quality based on the allowed number of reloads and intensity of demand for a location.

\subsection{Computational platform\label{sec:cp}}
The algorithm was implemented in Java, and CPLEX 12.6 was used to solve the linear RMP and integer pricing programs. Tree enumeration and column management were implemented using the Java collection library. All the computational runs were performed on an ACPI x64 machine (Intel Xeon E5630 processor @ 2.54 GHz, 12 GB RAM). The computation times reported are expressed in seconds, and we imposed a time limit of 7,200 seconds. The performance of the algorithm was tested on two different classes of test instances which are derived from real world data. The results from B\&P are benchmarked with the holistic EP model. The computational results for the EP model are reported.

\subsection{Instance generation\label{sec:ig}}
We generated two classes of test instances, A and B. The classes differ in the total number of destinations, \emph{i.e.}, dealers' locations. The solution for the loading problem depends on the number and types of vehicles, and the order of route sequence for a dealer location. Hence, we created two types of demand for dealer locations, high and low number of vehicles per location. Class A instances are low demand instances where $|S|\in \{15, 20, 25\}$, and class B instances are high demand instances with $|S| \in \{5,7,10\}$. In reality, class B represents city loads, few dealer locations with high demand for each of them, and class A represents non-city loads with less volume for high number of locations. For both the classes, the number of vehicles $|V|$ takes a value in the set $\{100, 200, 400 , 600\}$, and a total of $28$ different vehicle types were considered, \emph{i.e.}, the vehicles in the set $V$ consists of $28$ different vehicle types. Table \ref{tab:dims} represents the number of truck, sedan, and hatchback within each of the vehicle set. Accurate dimensions and weight of each vehicle type were obtained from OEMs, and other third party providers. Specifications for the auto-carrier like ramp lengths, maximum allowable slide angle, split ramps, and heights for each auto-carrier type were obtained from their computer-aided-design (CAD) drawings. Based on the available dimensions, the derived parameters used in LDP($t$) formulation are computed using trigonometry and force balance equations (see Appendix). The maximum number of allowable violations, ($\text{v}^t_{max}$) is a user parameter. We performed a computational study with $\text{v}^t_{max} \in \{0,2,4\}$. When $\text{v}^t_{max}$ takes a value $0$, a LIFO policy is imposed. We created $36$ class A and B instances.  For the results tables, the column headings are defined as follows, `Name' is the instance name, `$|V|$' is the number of vehicles, `$|S|$' is the number of dealers' locations, `$\text{v}^t_{max}$' is the maximum number of allowable violations, `opt' is the optimal objective value (total number of loads), `LR' is the load efficiency ratio, and `EP \%' is the final relative mixed-integer programming (MIP) gap reported by CPLEX from the holistic model EP1 after the stipulated runtime. Apart from PPL, another metric useful to examine the efficiency of a load is the ratio between the number of vehicles $|V|$ and number of auto-carriers used for loading (`opt' column) denoted by load efficiency ratio ($LR$). For the computational study, we used an auto-carrier type with nine ramps. The reason is to allow us to a make a comparison of performance of algorithm with respect to different classes of instances and EP1. However, as denoted by the formulation, the methodology can be used by auto-carrier with different capabilities and capacities. Furthermore, the vehicles for each dealer location are randomly assigned. For EP1 model, the number of available auto-carriers is given as $|V|/6$, i.e., a $LR$ of 6 is used. In the tables \ref{tab:resultsA} and \ref{tab:resultsB}, the higher MIP gap \% for the instances in classes A \& B using EP1 formulation indicates the necessity of a decomposition algorithm for ATP.

\begin{table}
	\centering
	\begin{tabular}{|c|c|c|c|}
		\hline	
		$|S|$ &  Truck  & Sedan & Hatchback \\
		\hline
		100 &  33  & 31 & 36 \\
		200 &  80  & 73 & 47 \\
		400 &  251  & 70 & 79 \\
		600 &  261  & 222 & 117 \\
		\hline
	\end{tabular}
	\caption{Distribution of Vehicle Types}
	\label{tab:dims}	
\end{table}

\subsection{Tests with Class A\label{sec:classa}}
Table \ref{tab:resultsA} summarizes the computational behavior of B\&P algorithm on class A instances. In this class the number of vehicles for each dealer location is less compared to the other class. For instance with a demand of less than or equal to 200 vehicles, the performance of the algorithm is very impressive. The holistic model was not able to generate a feasible solution for instances with 400 or more vehicles. Since the number of vehicles for a location is less (for A-100-15 with $|S|$ = 15, each location's demand is 6 to 7 vehicles on an average), B\&P algorithm uses the advantage of reloads well. The change in number of auto-carriers for instances having 600 vehicles with and without using reloads is around five to six auto-carriers. Also in general, with an increase in the percentage of trucks for an instance, $LR$ reduces as there are more number of exclusions during the formation of the set $\boldsymbol{\alpha}^t(R,J,V,L,M)$.

\begin{figure}
	\centering
	\includegraphics[scale=0.40]{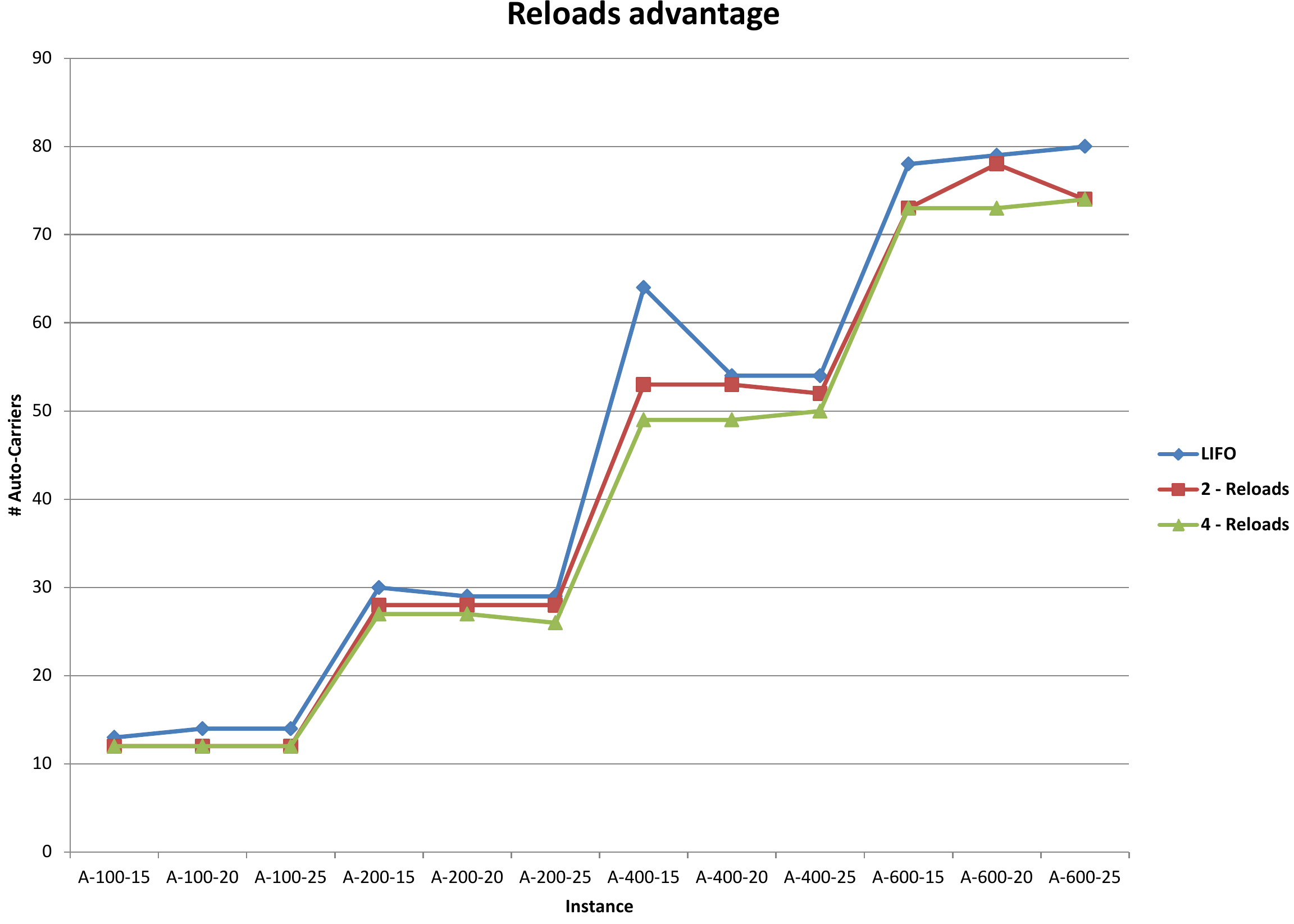}
	\caption{Advantage of reloads with class A instances}
	\label{fig:reloadA}
\end{figure}

Figure \ref{fig:reloadA} depicts the advantage of using reloads for class A instances, and the advantage increases with the increase in the number of vehicles for the instances. With more number of vehicles, B\&P approach had more opportunities to shuffle the vehicle-ramp assignment. Also, there is no pattern in terms of runtime with respect to a change in $\text{v}^t_{max}$, however runs with higher reloads were generally faster.

\subsection{Tests with Class B\label{sec:classb}}
Table \ref{tab:resultsB} summarizes the computational behavior of B\&P algorithm on class B instances. Only instance B-100-5 with LIFO was able to reach optimality with EP1 formulation. In class B, we noticed that the number of required auto-carriers did not significantly reduced with the increase in the number of allowed reloading. This is due to the reason that with a large volume of vehicles for each dealer's location, the loading problem was able to generate loads without utilizing the advantage of reloading. However, in class A due to the high number of dealers' locations for each instance, the number of vehicles for each location is low. Hence, the algorithm needed the allowance from reloading to reduce the total number of auto-carriers required for delivery. The high `$LR$' for class B also indicates the needlessness for a reload allowance as each location has a high demand for vehicles. This gives the opportunity for B\&P to build loads exclusively to a single dealer location. Also, class B instances had much better runtime than class A instances.

\begin{landscape}
	\begin{table}
		\begin{tabular}{cccccccc}
			\hline
			\noalign{\vskip0.15cm}
			Name & $|V|$ & $|S|$ & $\text{v}^t_{max}$ & Time & opt & $LR$ & EP\% \tabularnewline[0.15cm]
			\hline
			\noalign{\vskip\doublerulesep}
			\noalign{\vskip0.1cm}
			A-100-15 & 100 & 15 & 0 & 144 & 13 & 7.69 & 20.63\tabularnewline
			A-100-15 & 100 & 15 & 2 & 592 & 12 & 8.33 & 34.64\tabularnewline
			A-100-15 & 100 & 15 & 4 & 175 & 12 & 8.33 & 30.56\tabularnewline
			A-100-20 & 100 & 20 & 0 & 145 & 14 & 7.14 & 14.53\tabularnewline
			A-100-20 & 100 & 20 & 2 & 199 & 12 & 8.33 & 38.27\tabularnewline
			A-100-20 & 100 & 20 & 4 & 377 & 12 & 8.33 & 25.93\tabularnewline
			A-100-25 & 100 & 25 & 0 & 279 & 14 & 7.14 & 25.93\tabularnewline
			A-100-25 & 100 & 25 & 2 & 253 & 12 & 8.33 & 30.56\tabularnewline
			A-100-25 & 100 & 25 & 4 & 337 & 12 & 8.33 & 32.66\tabularnewline
			Avg   & & & & 277 & & 7.99 & 28.19\tabularnewline		
			\multicolumn{7}{c}{} \tabularnewline
			A-200-15 & 200 & 15 & 0 & 3,756 & 30 & 6.66 & 44.44\tabularnewline
			A-200-15 & 200 & 15 & 2 & 1,747 & 28 & 7.14 & 56.43\tabularnewline
			A-200-15 & 200 & 15 & 4 & 969 & 27   & 7.40 & --\tabularnewline
			A-200-20 & 200 & 20 & 0 & 1,903 & 29 & 6.89 & 45.80\tabularnewline
			A-200-20 & 200 & 20 & 2 & 1,790 & 28 & 7.14 & 50.06\tabularnewline
			A-200-20 & 200 & 20 & 4 & 896 & 27   & 7.40 & 51.69\tabularnewline
			A-200-25 & 200 & 25 & 0 & 1,274 & 29 & 6.89 & 50.62\tabularnewline
			A-200-25 & 200 & 25 & 2 & 1,976 & 28 & 7.14 & 56.43\tabularnewline
			A-200-25 & 200 & 25 & 4 & 1,127 & 26 & 7.69 & --\tabularnewline
			Avg   & & & & 1,715 & & 7.15 & 50.78\tabularnewline
			\hline \tabularnewline
		\end{tabular} \qquad 
		\begin{tabular}{ccccccccc}
			\hline
			\noalign{\vskip0.15cm}
			Name & $|V|$ & $|S|$ & $\text{v}^t_{max}$ & Time & opt & $LR$ & EP\%\tabularnewline[0.15cm]
			\hline
			\noalign{\vskip\doublerulesep}
			\noalign{\vskip0.1cm}
			A-400-15 & 400 & 15 & 0 & 1,872 & 64 & 6.25 & --\tabularnewline
			A-400-15 & 400 & 15 & 2 & 4,102 & 53 & 7.54 & --\tabularnewline
			A-400-15 & 400 & 15 & 4 & 1,288 & 49 & 8.16 & --\tabularnewline
			A-400-20 & 400 & 20 & 0 & 6,999 & 54 & 7.40 & --\tabularnewline
			A-400-20 & 400 & 20 & 2 & 7,200 & 53 & 7.54 & --\tabularnewline
			A-400-20 & 400 & 20 & 4 & 2,957 & 49 & 8.16 & --\tabularnewline
			A-400-25 & 400 & 25 & 0 & 3,692 & 54 & 7.40 & --\tabularnewline
			A-400-25 & 400 & 25 & 2 & 7,200 & 52 & 7.69 & --\tabularnewline
			A-400-25 & 400 & 25 & 4 & 4,581 & 50 & 8.00 & --\tabularnewline
			Avg & & & & 4,587 & & 7.57 & --\tabularnewline
			\multicolumn{6}{c}{} \tabularnewline
			A-600-15 & 600 & 15 & 0 & 6,230 & 78 & 7.69 & --\tabularnewline
			A-600-15 & 600 & 15 & 2 & 6,834 & 73 & 8.21 & --\tabularnewline
			A-600-15 & 600 & 15 & 4 & 5,878 & 73 & 8.21 & --\tabularnewline
			A-600-20 & 600 & 20 & 0 & 6,935 & 79 & 7.59 & --\tabularnewline
			A-600-20 & 600 & 20 & 2 & 6,633 & 78 & 7.69 & --\tabularnewline
			A-600-20 & 600 & 20 & 4 & 2,826 & 73 & 8.21 & --\tabularnewline
			A-600-25 & 600 & 25 & 0 & 7,063 & 80 & 7.50 & --\tabularnewline
			A-600-25 & 600 & 25 & 2 & 3,580 & 74 & 8.10 & --\tabularnewline
			A-600-25 & 600 & 25 & 4 & 2,905 & 74 & 8.10 & --\tabularnewline
			Avg & & & & 5,884 & & 7.92 & --\tabularnewline[0.1cm]
			\hline \tabularnewline
		\end{tabular}
		\caption{Computational results for Class A}
		\label{tab:resultsA}		
	\end{table}
\end{landscape}

\begin{landscape}
	\begin{table}
		\begin{tabular}{cccccccc}
			\hline
			\noalign{\vskip0.15cm}
			Name & $|V|$ & $|S|$ & $\text{v}^t_{max}$ & Time & opt & $LR$ & EP\% \tabularnewline[0.15cm]
			\hline
			\noalign{\vskip\doublerulesep}
			\noalign{\vskip0.1cm}
			B-100-5 & 100 & 5 & 0 & 42 & 12 & 8.33 & 0.00\tabularnewline
			B-100-5 & 100 & 5 & 2 & 41 & 12 & 8.33 & 39.94\tabularnewline
			B-100-5 & 100 & 5 & 4 & 46 & 12 & 8.33 & 23.37\tabularnewline
			B-100-7 & 100 & 7 & 0 & 41 & 13 & 7.69 & 14.53\tabularnewline
			B-100-7 & 100 & 7 & 2 & 47 & 12 & 8.33 & 14.53\tabularnewline
			B-100-7 & 100 & 7 & 4 & 94 & 12 & 8.33 & 34.64\tabularnewline
			B-100-10 & 100 & 10 & 0 & 226 & 13 & 7.69 & 14.53\tabularnewline
			B-100-10 & 100 & 10 & 2 & 72 & 13 & 7.69 & 20.63\tabularnewline
			B-100-10 & 100 & 10 & 4 & 70 & 13 & 7.69 & 25.93\tabularnewline
			Avg   &  & & & 75 & & 8.04 & 23.53\tabularnewline
			\multicolumn{6}{c}{} \tabularnewline
			B-200-5 & 200 & 5 & 0 & 190 & 24 & 8.33 & 38.27\tabularnewline
			B-200-5 & 200 & 5 & 2 & 284 & 24 & 8.33 & 52.72\tabularnewline
			B-200-5 & 200 & 5 & 4 & 273 & 24 & 8.33 & 38.27\tabularnewline
			B-200-7 & 200 & 7 & 0 & 282 & 24 & 8.33 & 36.51\tabularnewline
			B-200-7 & 200 & 7 & 2 & 377 & 24 & 8.33 & 50.06\tabularnewline
			B-200-7 & 200 & 7 & 4 & 712 & 24 & 8.33 & 54.18\tabularnewline
			B-200-10 & 200 & 10 & 0 & 663 & 24 & 8.33 & 32.66\tabularnewline
			B-200-10 & 200 & 10 & 2 & 456 & 24 & 8.33 & --\tabularnewline
			B-200-10 & 200 & 10 & 4 & 477 & 24 & 8.33 & 48.32\tabularnewline
			Avg   & & & & 412 & & 8.33 & 43.87\tabularnewline
			\hline \tabularnewline
		\end{tabular} \qquad 
		\begin{tabular}{cccccccc}
			\hline
			\noalign{\vskip0.15cm}
			Name & $|V|$ & $|S|$ & $\text{v}^t_{max}$ & Time & opt & $LR$ & EP\%\tabularnewline[0.15cm]
			\hline
			\noalign{\vskip\doublerulesep}
			\noalign{\vskip0.1cm}
			B-400-5 & 400 & 5 & 0 & 791 & 50 & 8.00  & --\tabularnewline
			B-400-5 & 400 & 5 & 2 & 804 & 49 & 8.16 & --\tabularnewline
			B-400-5 & 400 & 5 & 4 & 586 & 49 & 8.16 & --\tabularnewline
			B-400-7 & 400 & 7 & 0 & 703 & 49 & 8.16 & --\tabularnewline
			B-400-7 & 400 & 7 & 2 & 1,191 & 49 & 8.16 & --\tabularnewline
			B-400-7 & 400 & 7 & 4 & 1,464 & 49 & 8.16 & --\tabularnewline
			B-400-10 & 400 & 10 & 0 & 749 & 49 & 8.16 & --\tabularnewline
			B-400-10 & 400 & 10 & 2 & 1,429 & 49 & 8.16 & --\tabularnewline
			B-400-10 & 400 & 10 & 4 & 1,277 & 49 & 8.16 & --\tabularnewline
			Avg   & & & & 999 & & 8.14 & --\tabularnewline
			\multicolumn{6}{c}{} \tabularnewline
			B-600-5 & 600 & 5 & 0 & 1,852 & 74 & 8.11 & --\tabularnewline
			B-600-5 & 600 & 5 & 2 & 1,745 & 74 & 8.11 & --\tabularnewline
			B-600-5 & 600 & 5 & 4 & 2,392 & 73 & 8.21 & --\tabularnewline
			B-600-7 & 600 & 7 & 0 & 1,434 & 73 & 8.21 & --\tabularnewline
			B-600-7 & 600 & 7 & 2 & 2,693 & 73 & 8.21 & --\tabularnewline
			B-600-7 & 600 & 7 & 4 & 2,918 & 73 & 8.21 & --\tabularnewline
			B-600-10 & 600 & 10 & 0 & 1,835 & 74 & 8.11 & --\tabularnewline
			B-600-10 & 600 & 10 & 2 & 3,157 & 73 & 8.21 & --\tabularnewline
			B-600-10 & 600 & 10 & 4 & 2,716 & 74 & 8.11 & --\tabularnewline
			Avg   & & & & 2,304 & & 8.17 & --\tabularnewline[0.1cm]
			\hline \tabularnewline
		\end{tabular}
		\caption{Computational results for Class B}
		\label{tab:resultsB}
	\end{table}
\end{landscape}

\section{Conclusions}\label{sec:Conclusion}
Auto-carrier test transportation problem (ATP) addresses the challenges of optimal movement of vehicles from the auto-manufacturers to dealers' locations using auto-carriers. The two main components of ATP are routing and loading. Routing details the routes for the auto-carriers and loading suggests on how to load the individual vehicles in an auto-carrier. We present a heuristic for the routing problem, and depending on the routes generated, we provide an exact algorithm based on branch-and-price (B\&P) approach to solve auto-carrier loading problem. In most of the current practices for the loading problem, disjunctions for an assignment or numerical coefficients based on past experience are used for modeling purposes. In this work, we use actual dimensions of the vehicles for the loading problem, hence this work can be adopted for any generic vehicles without a reliance on the past data or previous experience to derive the numerical coefficients. Also, we relax the last-in-first-out (LIFO) policy during loading so the solution presents the trade-off in resource utilization and the number of reloads. The B\&P method for loading problem is tested with instances created from real-world data, and the efficiency of the model is evaluated with a holistic model using extensive computational experiments. For the loading problem, B\&P method outperforms holistic model in the computational experiments, and large scale instances are solved within two hours time limit. 

Future work include using routing problem along with loading in the B\&P framework. Future work include using routing problem along with loading in the B\&P framework. Multi-depot fuel constrained multiple vehicle routing problem introduced in \cite{sundar2015formulations} and \cite{sundar2016exact} could be appropriate. Another challenge faced by the logistics companies is capacity planning for auto-carriers, since auto-carriers are expensive assets. Given the stochastic nature of demand for the vehicles, the loading problem can be extended to a stochastic setting for capacity planning. From a computational perspective, the pricing problem can be solved independently so an implementation capable of solving pricing problems in parallel should also be evaluated.

\section*{Appendix}
\renewcommand{\thesection}{A}
\subsection{Calculation of effective weights at the steer, drive and tandem axle for a given configuration}
Given a configuration $\alpha \in \boldsymbol{\alpha}^t(R(t)\cup R^c_{\mathcal{SP}}(t),J,V,L,M)$ for an auto-carrier of type $t$, the effective weights acting at the steer, drive, and the tandem axles of the auto-carrier are calculated using force balance and moment balance equations. For the calculation of the effective weights at the three axles for a given configuration, the following assumptions are made.
\begin{enumerate}
	\item We ignore the slide angles $L$, \emph{i.e.}, we calculate the weights assuming all the ramps are perfectly horizontal. 
	\item When a vehicle is loaded on a ramp, the force due to the entire weight of the vehicle acts at a point on the ground that is vertically below the geometric center of the ramp. 
	\item In case of a split ramp, \emph{i.e.}, one vehicle loaded on two ramps, the force acts at a point on the ground vertically below the geometric center of the combined ramp. 
	
\end{enumerate}
\begin{figure}
	\centering
	\includegraphics[scale=1]{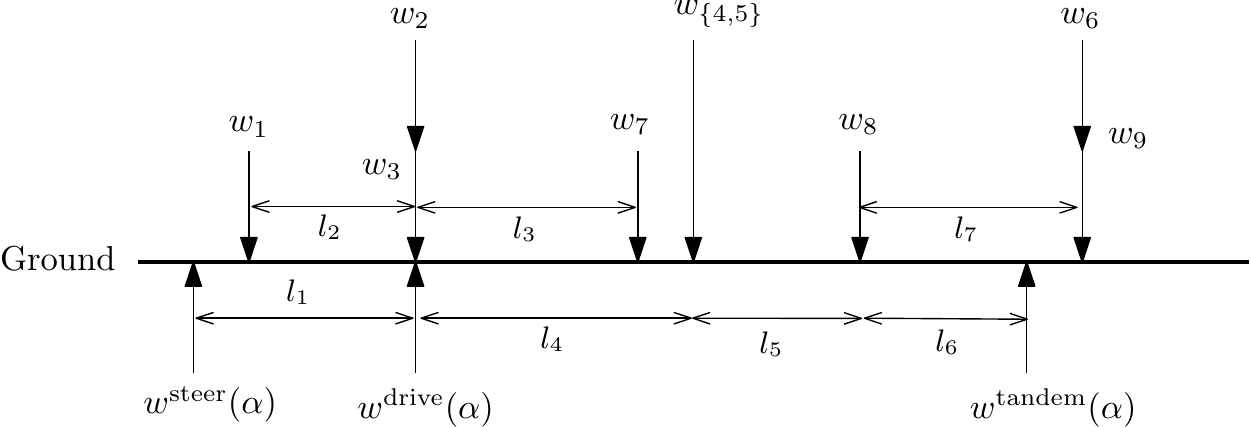}
	\caption{Free-body diagram for an auto-carrier.}
	\label{fig:weight}
\end{figure}
The lengths, geometric center locations, and the distance from geometric center for the ramps to every axle are obtained from the auto-carrier's CAD drawings. Figure \ref{fig:weight} shows the free-body diagram. The force due to the weight of each vehicle in a given configuration acts at points which are vertically below their corresponding geometric center on the ground. Let $w_i$ denotes the force due to the weight of the vehicle loaded on the ramp $i$ (note that $i$ can also denote a split ramp, as is the case of $w_{\{4,5\}}$). The unknown reaction forces at steer, drive, and tandem axles are indicated by $w^{\text{steer}}(\alpha)$, $w^{\text{drive}}(\alpha)$, and $w^{\text{tandem}}(\alpha)$, respectively. The lengths $l_k$ ($k=1,\dots,7$) are known parameters. The force balance equation for the set of forces is given by \[ w^{\text{steer}}(\alpha) + w^{\text{drive}}(\alpha) + w^{\text{tandem}}(\alpha) = w_1 + w_2 + w_3 + w_{\{4,5\}} + w_6 + w_7 + w_8 + w_9. \] The moment balance principle states that the sum of the clockwise moments about a given point is equal to the sum of the anti-clockwise moments about the same point. The moment balance equation about the point at which the reaction force $w^{\text{drive}}(\alpha)$ acts is given by
\[l_1 \cdot w^{\text{steer}}(\alpha) + l_3 \cdot w_7 + l_4 \cdot w_{\{4,5\}} + (l_4+l_5)\cdot w_8 + (l_4+l_5+l_7)\cdot (w_6+w_9) = l_1\cdot w_1 + (l_4 + l_5 + l_6)\cdot w^{\text{tandem}}(\alpha). \] Similarly, moment balance equations for $w^{\text{steer}}(\alpha)$ and $w^{\text{tandem}}(\alpha)$ can be derived. The system of equations are then solved to compute the unknown reaction forces, \emph{i.e.}, $w^{\text{steer}}(\alpha)$, $w^{\text{drive}}(\alpha)$, and $w^{\text{tandem}}(\alpha)$.

\subsection{Calculation of height parameters for a given configuration in an auto-carrier}
\begin{figure}
	\centering
	\includegraphics[scale=0.8]{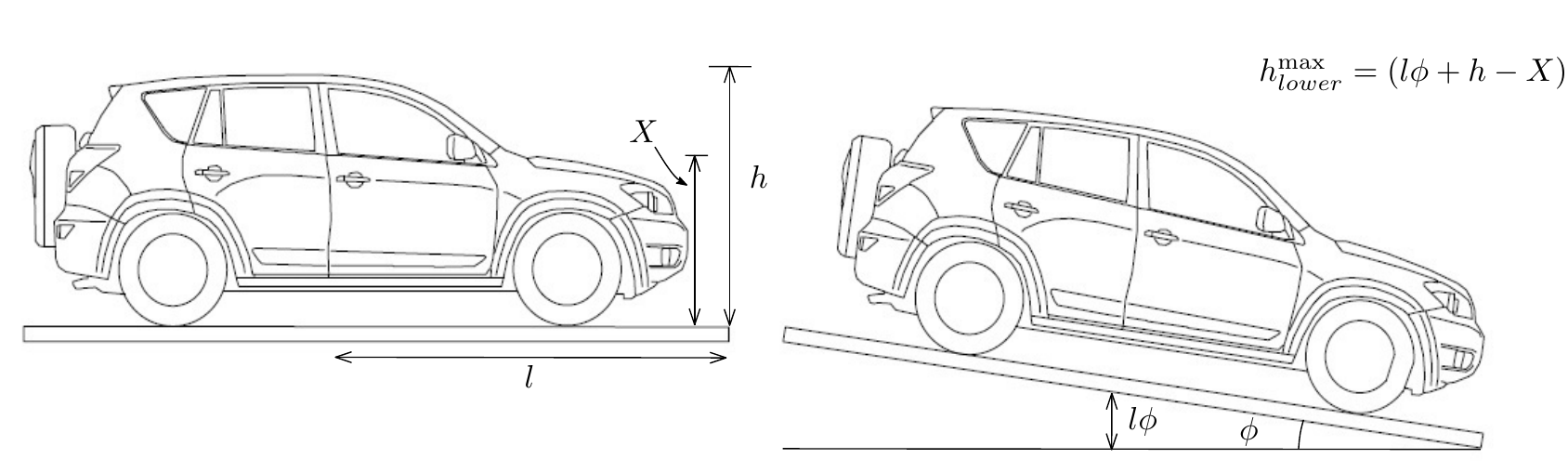}
	\caption{Height parameter calculation for a given configuration.}
	\label{fig:height}
\end{figure}

Given a configuration $\alpha \in \boldsymbol{\alpha}^t(R,J,V,L,M)$ where $R \in R_{\mathcal{U}}(t)$, \emph{i.e.}, a configuration for a ramp in the upper deck, we detail the procedure to calculate $h_{\text{upper}}^{\text{gain}}(\alpha)$. Any ramp in the auto-carrier of type $t$ can slide within a maximum allowable slide angle about its pivot. A ramp can have pivots at both its ends. The configuration $\alpha$ specifies a slide angle, a loaded vehicle and the direction of slide. The direction of slide  specifies the pivot of rotation. For instance, consider the ramp shown in Figure \ref{fig:height} where $h$ is the maximum height of the loaded vehicle on the ramp, and $l$ is the distance from the pivot to the point on the ramp where the vehicle attains this maximum height. Both these values can be computed with the knowledge of the vehicle and ramp dimensions. Now, if the slide angle is $\phi$, then the gain in height  produced by sliding the ramp is computed approximately as $l \sin \phi \approx l \phi$ for small angles $\phi$. Hence $h_{\text{upper}}^{\text{gain}}(\alpha) \approx l \phi.$

Now, we detail the procedure to compute $h_{\text{lower}}^{\text{max}}(\alpha)$ for a configuration $\alpha \in \boldsymbol{\alpha}^t(R,J,V,L,M)$ where $R \in  R_{\mathcal{L}}(t)$. Suppose $i$ and $\phi$ denote the ramp and its angle of slide in the configuration $\alpha$, respectively, then the parameter $h_{\text{lower}}^{\text{max}}(\alpha)$ is the maximum slide for the ramp vertically above the ramp $i$, based on the vehicle in the ramp $i$ and the slide angle $\ell$. Then, $h_{\text{lower}}^{\text{max}}(\alpha) = (l \phi + h - X)$ as depicted in Figure \ref{fig:height}. 
\subsection{Selection of dealers' locations for route sequence}
Routing heuristic helps the planner to choose the set of dealers' locations near the destination $d$. The planner chooses the destination dealer $d$, viewing angle $\varphi$, and offset distance $\nu$, based on which a polygon is constructed, and then Lin-Kernighan heuristic is used to determine the route sequence for the vehicles to be delivered.

\begin{figure}[h]
	\centering
	\includegraphics[scale=0.3]{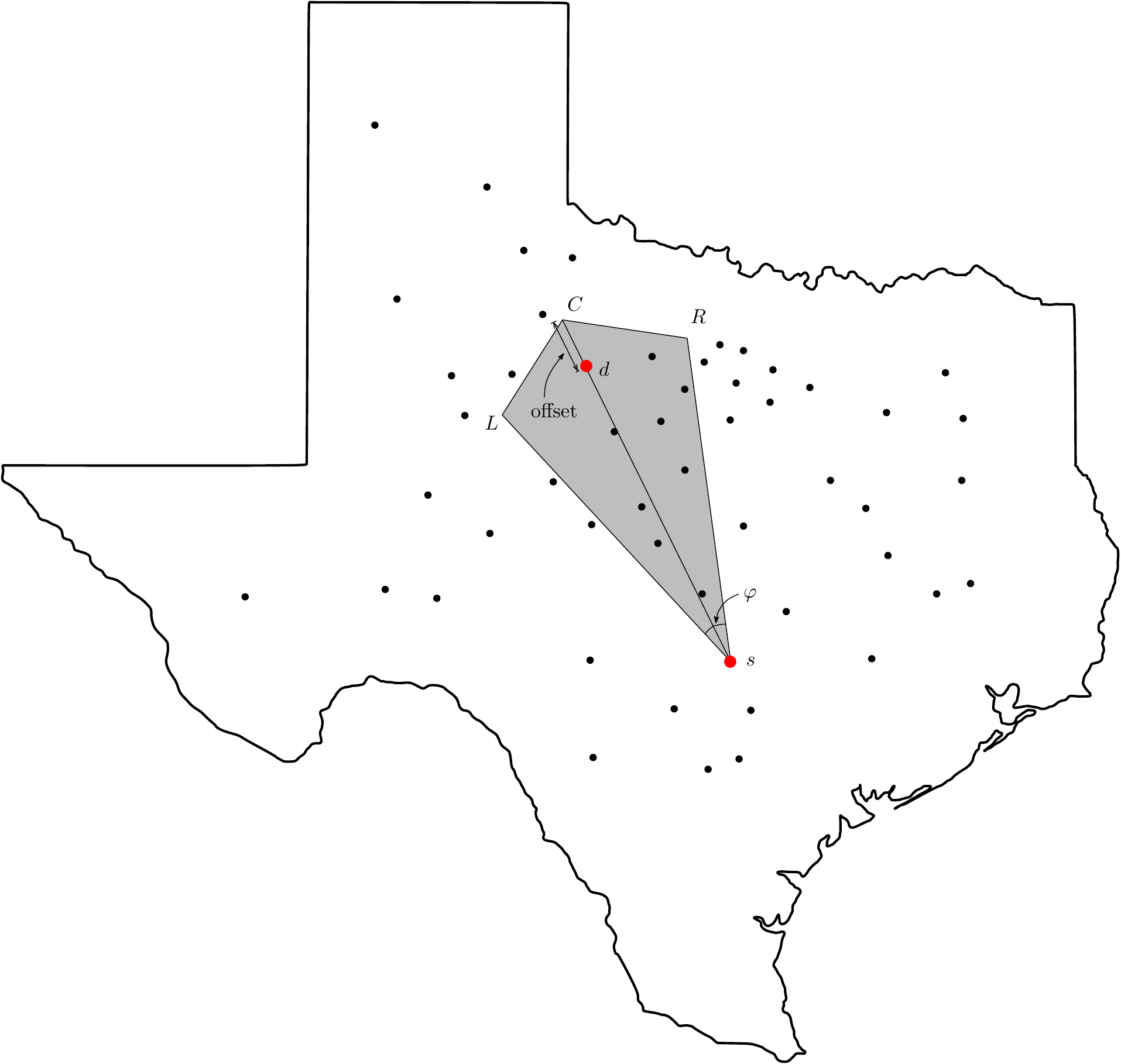}
	\caption{Illustration of Routing heuristics.}
	\label{fig:map}
\end{figure}

Let $s$ be the ALC's distribution center, $d$ be the target dealer location, $\varphi$ is the viewing angle, `$lat$' and `$lon$' represent latitude and longitude in radian unit of measure for a location, respectively. Also, $\sin$, $\cos$, $\rm atan2$, $\rm asin$ represent sine, cosine, arctangent, and inverse sine trigonometric functions, respectively. Let $\bar{y}$ be given as \[ \bar{y} = \sin(d_{lon}-s_{lon}).\cos(d_{lat}), \] and $\bar{x}$ is given as, \[ \bar{x} = \cos(s_{lat}).\sin(d_{lat})- \sin(s_{lat}).\cos(d_{lat}).\cos(d_{lon}-s_{lon}). \] Then, the bearing angles $\theta$ is given as, $\theta = \rm atan2(y,x)$. The distance $di$ between $s$ and $d$ is calculated using the Haversine formula given as follows, \[ a = \sin((d_{lat}-s_{lat})/2).\sin((d_{lat}-s_{lat})/2) + \sin((d_{lon}-s_{lon})/2).\sin((d_{lon}-s_{lon})/2)*\cos(s_{lat})*\cos(d_{lat}), \]
\[ c = 2.\rm atan2(\sqrt{a},\sqrt{1-a}). \] Then, the distance $di$ is given as $di = R.c$, where $R$ is the diameter of earth, a value of 6,371 kilometers is used. Using the values $di$ and $\theta$, the latitude and longitude for the three extreme points of polygon $\chi$ are calculated. Let $C$, $L$, and $R$ be the center, left, and right extreme points of the polygon $\chi$, respectively. The fourth extreme point is the distribution center $s$. The angle values of $\theta$, $\theta + \varphi /2$, and $\theta - \varphi /2$ represent the bearing angle for $C$, $L$, and $R$, respectively. Then the latitude and longitude for $C$ is calculated as \[ C_{lat} = {\rm{asin}} (\sin(s_{lat}.\cos(K/R) + \cos(s_{lat}).\sin(K/R).\cos(\theta)) \] and \[ C_{lon} = s_{lon} + {\rm{atan2}}(\sin(\theta).\sin(K/R).\cos(s_{lat}),\cos(K/R)-\sin(s_{lat}).\sin(C_{lat})), \] where $K=(di + \nu)/R$. Similarly, $L_{lat}$ and $L_{lon}$ can be determined using $\theta + \varphi /2$ instead of $\theta$, $R_{lat}$ and $R_{lon}$ can be determined using $\theta - \varphi /2$ instead of $\theta$ in the above formula. For calculating $L_{lat}$, $L_{lon}$, $R_{lat}$, and $R_{lon}$, we use $K=di/R$. Once we construct a polygon $\chi$ using the latitudes and longitudes of $C$, $L$, $R$, and $s$, the set of dealers' locations $S$ within the polygon $\chi$ is determined. The details of Haversine formula are given in \cite{Rick1999}.

%
\bibliography{nvl}
\end{document}